\documentclass[reqno]{amsart}

    \usepackage{float}
    \usepackage{amsmath}
    \usepackage{graphicx}
    \usepackage{latexsym}
    \usepackage{amsfonts}
    \usepackage{amssymb}
    \usepackage{verbatim}
    \usepackage{xcolor}
    \usepackage{url}
    \usepackage{hyperref}

    \theoremstyle{plain}
    \newtheorem{theorem}{Theorem}
    \newtheorem{corollary}[theorem]{Corollary}
    \newtheorem{lemma}[theorem]{Lemma}
    \newtheorem{proposition}[theorem]{Proposition}
    \theoremstyle{definition}
    \newtheorem{assumption}{Assumption}

    \newtheorem*{remark*}{Remark}

    \newcommand{\pr}{\mathbf P}
    \newcommand{\e}{\mathbf E}

\begin{document}
    \title[LOCAL LIMIT THEOREM FOR BRANCHING PROCESSES WITH IMMIGRATION ]
    {LOCAL LIMIT THEOREM FOR CRITICAL BRANCHING PROCESSES WITH IMMIGRATION}
    \author[]{Sadillo Sharipov, Alexander Tarasov and Vitali Wachtel}
    \address{V.I.Romanovskiy Institute of Mathematics, Tashkent, Uzbekistan}
    \email{sadi.sharipov@yahoo.com}
    \address{Hildesheim University, Germany}
    \email{tarasov@uni-hildesheim.de}
    \address{Faculty of Mathematics, Bielefeld University, Germany}
    \email{wachtel@math.uni-bielefeld.de}

\begin{abstract}
Let $\{Y_{n}, n\geq 0\}$ be a critical branching process with immigration. Under the minimal moment restrictions, namely, assuming finite variance for the number of particles of the offspring and finite mean for the number of immigrating particles, we obtain a local limit theorem for the process $\{Y_{n}\}$.
\end{abstract}


    \maketitle

\section{Introduction, main results and discussion}

Let $\left\{Y_{n}\right\}$ be a discrete-time branching process with the offspring distribution $\left\{p_{k} \right\}$
and with the immigration governed by the distribution $\left\{q_{k} \right\}$.
A mathematically rigorous description can be as follows. Let $\left\{\xi_{n,i}\right\}$ be independent, identically distributed nonnegative and integer valued random variables with the distribution $\left\{p_{k} \right\}$. Let also $\left\{\eta_n\right\}$ be independent, identically distributed nonnegative and integer valued random variables with the distribution $\left\{q_{k} \right\}$. Furthermore, we shall assume that $\left\{\xi_{n,i}\right\}$ and $\left\{\eta_n\right\}$ are
independent from each other. Then the process $\left\{Y_{n}\right\}$ is given by the following recursive equation
$$
Y_{n+1}=\sum_{j=1}^{Y_{n}} \xi_{n+1,j}+\eta_{n+1}, \quad n \geq 0.
$$
Here, $Y_{0}$ is either deterministic or random; in the latter case $Y_{0}$ is assumed independent of $\left\{\xi_{n,j}\right\}$ and
$\left\{\eta_n\right\}$.

We shall consider critical processes, that is,
$$ \mathbf{E}\xi_{1,1}=1. $$

In the case $q_0=1$, $\{Y_{n}\}$ reduces to a Galton-Watson process, which will be denoted by
$\{Z_{n}\}$. 
Local probabilities for critical branching processes have been extensively studied in the literature. Apparently, one of the first investigations of local probabilities was conducted by Zolotarev \cite{Zolotarev58}, who derived asymptotics, as time tends to infinity, for the probability that a critical Markov branching process in continuous time has a fixed size $k$. Chistyakov~\cite{Chistyakov58} studied continuous-time Markov branching processes and established local limit theorems for both critical and supercritical branching processes, assuming the fourth moment of the offspring distribution is finite.
Later, under assumption that $\mathbf{E}\left(\xi_{1,1}\right)^{4} < \infty$, Nagaev and Mukhamedkhanova \cite{NM66} derived a local limit theorem for a critical Galton-Watson process.
In Athreya and Ney's book \cite{AN72}, a local limit theorem was established for a critical Galton-Watson process under the hypothesis $\mathbf{E}\left(\xi_{1,1}\right)^{3+\delta} < \infty$ for some $\delta > 0$.

Concurrently, Kesten, Ney, and Spitzer \cite{KestenNeySpitzer67} obtained a local limit theorem for a critical Galton-Watson process. Specifically, they imposed the conditions that $\left\{p_{k}\right\}$ is aperiodic and
\begin{equation}\label{KNScondition}
\mathbf{E}\xi_{1,1}^{2} \log\left(1 + \xi_{1,1}\right) < \infty
\end{equation}
and showed that if $k$ and $n$ tend to infinity in such a way that the ratio $\frac{k}{n}$ remains bounded, then
\begin{equation}\label{KNSresult}
\lim_{n\to \infty} n^{2} \exp\left\{\frac{2k}{Bn}\right\} \pr\left(Z_{n} = k\right) = \frac{4}{B^{2}},
\end{equation}
where 
$$
B:=\mathbf{E}\xi_{1,1}(\xi_{1,1}-1).
$$

The authors of \cite{KestenNeySpitzer67} indicated that the relation \eqref{KNSresult} holds without excessive moment restrictions, i.e., it is sufficient to satisfy the condition $B < \infty$. However, their proof of \eqref{KNSresult} was carried out under the condition \eqref{KNScondition}. Nagaev and Vakhtel \cite{NV06} have later shown that the local limit theorem
for critical Galton-Watson processes is valid under the hypothesis that $B < \infty$. More precisely, they proved that
$$
\lim_{n\to \infty} \frac{B^{2}n^{2}}{4} \left(1 + \frac{2}{Bn}\right)^{k+1} \pr\left(Z_{n} = k\right) = 1
$$
for $k$ and $n$ such that $k/n$ remains bounded.

Let us now turn to critical processes with immigration. If the immigration is non-degenerate, then the state $0$ is not absorbing. If one additionally assumes that $B$ is finite, then the chain $|{Y_n|}$ is either null-recurrent or transient. Furthermore, the sequence
$\frac{2Y_n}{Bn}$ converges weakly towards the Gamma-distribution with parameters $\gamma=\frac{2\lambda}{B}$ and $1$, where $\lambda=\mathbf{E}\eta_1$. This convergence has been proven by
Seneta~\cite{Seneta70}.

Mellein \cite{Mellein82} has proved the corresponding local limit theorem. More precisely, imposing the restrictions
\begin{equation}\label{Mcondition}
\mathbf{E}\xi_{1,1}^{2}\log\left(1+\xi_{1,1}\right)<\infty \quad\text{and}\quad \mathbf{E}\eta_{1}\log\left(1+\eta_{1}\right)<\infty,
\end{equation}
he has shown that
\begin{equation}
\label{eq:loc-limit}
n\,\pr\left(Y_{n}=k\right) \sim
\frac{1}{\Gamma\left(\gamma\right)}\left(\frac{2k}{Bn}\right)^{\gamma-1}e^{-2k/(Bn)}
\end{equation}
as $n,k\to\infty$ and $\frac{k}{n}$ remains bounded. (Here and in what follows, the symbol $\sim$ denotes the asymptotic equivalence, that is, $a_n\sim b_n$ means that $\frac{a_n}{b_n}\to1$ as $n\to\infty.$)

Rahimov~\cite{Rahimov88} has proved an analogue of \eqref{eq:loc-limit} for critical branching processes with decreasing immigration. In his model, the random variables $\eta_n$ are not identically distributed and $\e \eta_n\to0$ as $n\to\infty$. Similar to \cite{KestenNeySpitzer67} and \cite{Mellein82}, it is assumed in~\cite{Rahimov88} that the offspring law satisfies \eqref{Mcondition}. Moreover, the author uses the assumption that the sequence $\e \eta^2_n$ is bounded.

The first moment condition in \eqref{Mcondition} is the same as in~\cite{KestenNeySpitzer67}.
Since the moment characteristics in \eqref{Mcondition} do not appear in \eqref{eq:loc-limit}, it is natural to expect
that \eqref{eq:loc-limit} should be valid without assumption \eqref{Mcondition}. This conjecture is supported also by the analogy with Galton-Watson processes.
But it turns out, see Sharipov and Wachtel \cite{ShW}, that this conjecture is wrong. The main result of that paper states that
\begin{equation}\label{ShW}
\pr\left(Y_{n}= k\right)\sim
\frac{1}{\Gamma(\gamma+1)}
\left(\frac{2}{B}\right)^{\gamma}
\frac{k^{\gamma-1} L\left(k\right)}{n^\gamma L\left(n\right)}
\end{equation}
as $n,k\to\infty$ and $k=o\left(n\right)$, where $L\left(x\right)$ is a slowly varying function.
$L\left(x\right)$ is asymptotically constant provided that the moment conditions in \eqref{Mcondition} hold.
Moreover, it may happen that $L(x)$ converges to zero or to infinity.

Since $L\left(x\right)$ is slowly varying, there exists a sequence $m_{n}=o\left(n\right)$ such that $L(m_{n})/L(n) \to 1$.
Thus, we infer from \eqref{ShW} that
$$
\pr\left(Y_{n}= k\right)\sim
\frac{1}{\Gamma(\gamma+1)}
\left(\frac{2}{B}\right)^{\gamma}
\frac{k^{\gamma-1}}{n^{\gamma}}
$$
provided that $k=o\left(n\right)$ and $k\geq m_{n}$.

This leads to the following conjecture: for every $\varepsilon>0$ one has
\begin{equation}\label{conjecture}
\sup_{k\geq \varepsilon n }\left| n\pr\left(Y_{n}=k\right)-\frac{1}{\Gamma\left(\gamma\right)}\frac{2}{B}\left(\frac{2k}{Bn} \right)^{\gamma-1}e^{-2k/(Bn)} \right| \to 0
\end{equation}
for \emph{all} critical processes with immigration, that is, \eqref{Mcondition} is not needed for the validity of the local limit theorem in the domain of normal deviations.

The main purpose of this paper is to proof the validity of this conjuncture.
\begin{theorem}\label{thm:loc.lim}
Assume that $\{p_{k}\}$ and $\{q_{k}\}$ are aperiodic. If $B<\infty$ and $\lambda<\infty$ then there exists a sequence $\varepsilon_n\to0$ such that 
\begin{equation}\label{conjecture}
\sup_{k\geq n\varepsilon_n }\left| n\pr\left(Y_{n}=k\right)-\frac{1}{\Gamma\left(\gamma\right)}\frac{2}{B}\left(\frac{2k}{Bn} \right)^{\gamma-1}e^{-2k/(Bn)} \right| \to 0
\quad\text{as }n\to\infty.
\end{equation}
\end{theorem}
As it has been proved in \cite{ShW}, see Remark 1 there, Theorem~\ref{thm:loc.lim}
implies that in the case $\gamma>1$ one has the most standard form of the local limit theorem.
\begin{corollary}
Assume that the assumptions of Theorem~\ref{thm:loc.lim} hold and that $\gamma>1$. Then
\begin{equation*}
\sup_{k\geq 0 }\left| n\pr\left(Y_{n}=k\right)-\frac{1}{\Gamma\left(\gamma\right)}\frac{2}{B}\left(\frac{2k}{Bn} \right)^{\gamma-1}e^{-2k/(Bn)} \right| \to 0
\quad\text{as }n\to\infty.
\end{equation*}
\end{corollary}
Let us explain our motivation for studying local probabilities for branching processes with immigration. 

First, it is rather natural to clarify the role of the extra moment conditions in~\eqref{Mcondition}, which have been introduced in Mellein's paper~\cite{Mellein82}. It was shown in~\cite{ShW} that these restrictions are needed to describe asymptotic properties of local lower deviation probabilities $\pr(Y_n=k)$ for $k=o(n)$. In particular, it has been proved there that if $\eqref{Mcondition}$ is violated then the local probabilities may significantly deviate from the density of the limiting law. Our Theorem~\ref{thm:loc.lim} shows that the moment restrictions in \eqref{Mcondition} are not needed in the domain of normal deviations, that is, for $k$ comparable with $n$.

Second, critical branching processes are closely related to Markov chains with asymptotically zero drift. Applying, for example, the results of Section~11.3 from~\cite{DKW-book}, we infer that the chain 
$$
X_n:=\sqrt{Y_n},\quad n\ge 0
$$
satisfies, as $x\to\infty$,
\begin{align}
\label{eq:lamperti}
\e[X_1-X_0|X_0=x]\sim\frac{\mu}{x}
\quad\text{and}\quad
\e[(X_1-X_0)^2|X_0=x]\sim b,
\end{align}
where
$$
\mu=\frac{\lambda-B/4}{2}
\quad\text{and}\quad 
b=\frac{B}{4}.
$$
Furthermore, the distribution of $X_1-X_0$ conditioned on $X_0=x$ converges, as $x\to\infty$, towards the normal law with mean $0$ and variance $B/4$.

It is a very hard task to prove local limit theorems for Markov chains satisfying \eqref{eq:lamperti}. Only some particular cases have been studied in the literature. Rosenkrantz~\cite{Rosenkrantz66} has proved a local limit theorem for a chain on $\{0,1,2,\ldots\}$ with nearest neighbors jumps and with the transition probabilities 
\begin{align*}
&\pr(X_1=x+1|X_0=x)=\frac{1}{2}\left(1+\frac{\lambda}{x+\lambda}\right),\\
&\pr(X_1=x+1|X_0=x)=\frac{1}{2}\left(1-\frac{\lambda}{x+\lambda}\right).
\end{align*}
The proof is based on the representation of local probabilities in terms of orthogonal polynomials.
This approach can be extended to Markov chains on one-dimensional hypergroups, see Ehring~\cite{Ehring94} and Guillotin-Plantard~\cite{Plantard11}. In all these cases, the connection to orthogonal polynomials allows one to use standard Fourier methods in the proof of local limit theorem. But it is absolutely unclear how to use Fourier transforms in a general situation.

Alexander~\cite{Alex11} has derived a local limit theorem for null-recurrent chains with nearest neighbors jumps. His proof is based on embeddings of a discrete time chain into a bessel process.

An alternative method for proving local limit theorems for Markov chains has been developed by Konakov and Mammen, see \cite{KM00} and \cite{KM01}. They have shown that if a properly rescaled Markov chain converges to a diffusion with smooth drift and diffusion coefficients, then the corresponding local limit theorem is valid. If a non-positive (transient or null-recurrent) chain $\{X_n\}$  satisfies \eqref{eq:lamperti} then $\frac{X_{nt}}{\sqrt{n}}$ converges typically to a Bessel process, see Section~4.7 in \cite{DKW-book}. Since the drift of a Bessel process is proportional to $\frac{1}{x}$, the method of Konakov and Mammen does not work in this case.
(Note that the singularity of the drift can formally be removed by considering the squared chain which will converge to a squared Bessel process. But it does not solve the problem, since the diffusion coefficient of every squared Bessel process is proportional to $\sqrt{x}$ and this function is not smooth at zero.) So, one needs a completely new approach to local asymptotics for Markov chains with asymptotically zero drift. Here we come to our main motivation which was to test the following idea. Fix some $\varepsilon\in(0,1)$ and set $m=\lfloor(1-\varepsilon^2)n\rfloor$. By the Markov property at time $m$,
\begin{align*}
\pr(X_n=x)
&=\sum_y\pr(X_m=y)\pr(X_{n-m}=x|X_0=y)\\
&=\sum_{y:|x-y|\le \varepsilon^{1/2}\sqrt{n}}\pr(X_m=y)\pr(X_{n-m}=x|X_0=y)\\
&\hspace{2cm}+\sum_{y:|x-y|> \varepsilon^{1/2}\sqrt{n}}\pr(X_m=y)\pr(X_{n-m}=x|X_0=y).
\end{align*}
If $\frac{X_{nt}}{\sqrt{n}}$ converges to a Bessel process then the time $n-m\sim \varepsilon^2 n$ is too small to get from $y$ to $x$ for $|x-y|> \varepsilon^{1/2}\sqrt{n}$. This leads naturally to the relation 
\begin{align*}
\pr(X_n=x)\approx \sum_{y:|x-y|\le \varepsilon^{1/2}\sqrt{n}}\pr(X_m=y)\pr(X_{n-m}=x|X_0=y).
\end{align*}
If we additionally assume that $\pr(X_1-X_0\in\cdot|X_0=x)$ converges, as $x\to\infty$, to a proper distribution with zero mean and finite variance, then one can expect that the behaviour of $\pr(X_{n-m}=x|X_0=y)$ can be described by the standard local \emph{central limit theorem}:
\begin{align*}
\pr(X_{n-m}=x|X_0=y)\approx \frac{1}{\sqrt{2\pi b(n-m)}}\exp\left\{-\frac{(x-y)^2}{2b(n-m)}\right\}.
\end{align*}
Therefore,
$$
\sqrt{n}\pr(X_n=x)\approx \frac{1}{\varepsilon\sqrt{2\pi b}}
\e\left[\exp\left\{-\frac{(x-X_{m})^2}{2b\varepsilon^2 n}\right\}\right].
$$
Using first the limit theorem for $\{X_n\}$ and letting then $\varepsilon\to0$, one concludes that the right hand sided of this relation is asymptotically equivalent to the values of the limiting density at the point $\frac{x}{\sqrt{bn}}$.

To implement this idea rigorously one needs two ingredients
\begin{itemize}
    \item[(a)] appropriate upper bounds for local large deviation probabilities\\ $\pr(X_{n-m}=x|X_0=y)$ with
    $|x-y|> \varepsilon^{1/2}\sqrt{n}$;
    \item[(b)] a local central limit theorem for $\pr(X_{n-m}=x|X_0=y)$\\ with 
    $|x-y|\le \varepsilon^{1/2}\sqrt{n}$.
\end{itemize}
It is not clear at the moment, how to deal with these problems in a general case.
But the proof of our Theorem~\ref{thm:loc.lim}, which implements the above strategy in the case of branching processes, confirms the capability of the approach in the case of Markov chains with sufficiently rich internal structure. 

A similar strategy has been used in \cite{DW15} to prove local limit theorems for conditioned multidimensional random walks. In this case one could use the standard bounds for concentration functions of random walks and local limit theorems for unconditioned walks with i.i.d. increments.

The remaining part of the paper is organized as follows. In Section 2 we establish some estimates for the concentration function of processes with immigration. In Section 3, we approximate Galton-Watson process with a branching process which offspring follows geometric law. Section 4 contains some results on upper bound for the local probabilities of Galton-Watson processes. Finally, in Section 5, we provide the proof of the local limit theorem. Throughout the proofs, we denote by $C$, $C_0$, $C_1,\ldots$ the generic positive numbers which may change their values from place to place and may depend only on distributions $\{p_k\}$ and $\{q_k\}$.  

In the proofs we shall always (without mentioning it) assume that the distributions $\{p_k\}$ and $\{q_k\}$ are aperiodic and that $\lambda,B\in(0,\infty).$

\section{Estimates for the concentration function of $Y_n$}
In this section, we derive several estimates for mass functions of branching processes with and without immigration. These estimates will play a crucial role in the proof of our local limit theorem.

\subsection{Local estimates for critical Galton-Watson processes.}
Let $\{Z_n\}$ be a critical Galton-Watson process with the offspring generating function 
$$
f(s):=\e s^{\xi_{1,1}}=\sum_{k=0}^\infty p_ks^k,\quad |s|\le 1.
$$
\begin{lemma}\label{lem:max-bound}
If $m=o(k)$, then
$$
k\max_{1\leq r\leq k, j\geq ck }\pr(Z_{m}=r | Z_{0}=j) \to 0
\quad\text{as }m,k\to\infty
$$
for any $c>2(1-e^{-1/2})^{-1}$. 

Moreover, there exists a constant $C$ such that
$$
k\max_{1\leq r\leq k, j\geq ck }\pr(Z_{m}=r | Z_{0}=j) \le C
$$
for all $m\le k$ and any $c>2(1-e^{-1/2})^{-1}$.
\end{lemma}
\begin{proof}
Let $\{\xi_i^{(m)}\}_{i\geq 1}$ be a sequence of independent identically
distributed random variables whose common distribution is defined by
$$
\pr(\xi_{1}^{(m)}=r)
=
\pr(Z_{m}=r|Z_{0}=1, Z_{m}>0),
\qquad r\geq 1.
$$
Set
$
S_{j}^{(m)}=\sum_{i=1}^j \xi_i^{(m)}$.
By the total probability law,
\begin{align}
\label{eq:binom}
\pr(Z_{m}=r | Z_{0}=j)=\sum_{l=1}^{j} \binom{j}{l}Q_{m}^{l}(1-Q_{m})^{j-l}\pr(S_{l}^{(m)}=r),
\end{align}
where
$$
Q_m:=\pr(Z_m>0|Z_0=1).
$$
Let us now perform the exponential change of measure with parameter $(-Q_{m})$. This means that we introduce a new measure, which is given by the equality 
\begin{equation}\label{eq:exponentialchange}
 \widehat{\pr}(\xi_{1}^{(m)}=r)=\frac{e^{-Q_{m}r}}{\mathbf{E}[e^{-Q_{m}\xi_{1}^{(m)}}]}\pr(\xi_{1}^{(m)}=r), \ \ r \geq 1. 
\end{equation}
Then
\begin{align}
\label{eq:exp.change}
\pr(S_{l}^{(m)}=r)=e^{Q_{m}r}
\left(\mathbf{E}[e^{-Q_{m}\xi_{1}^{(m)}}]\right)^{l}\,\widehat{\pr}(S_{l}^{(m)}=r), \ \ r \geq 1.
\end{align}
Notice next that, by the definition of $\xi_{1}^{(m)}$,
$$
\mathbf{E}[e^{-Q_{m}\xi_{1}^{(m)}}]=\mathbf{E}[e^{-Q_{m}Z_{m}}| Z_{0}=1, Z_{m}>0].
$$
Applying the Yaglom limit theorem, we conclude that
$$
\mathbf{E}[e^{-Q_{m}\xi_{1}^{(m)}}] \to \frac{1}{2}, \ \ m\to \infty.
$$
Thus, for all $m$ large enough, say $m\geq m_{0}$,
\begin{equation}\label{eq:inequalityforexponetial}
e^{-1}\leq \mathbf{E}[e^{-Q_{m}\xi_{1}^{(m)}}]\leq e^{-1/2}.
\end{equation}
Applying first \eqref{eq:exp.change} and then \eqref{eq:inequalityforexponetial} to the right hand side of \eqref{eq:binom}, we obtain
\begin{align}\label{mainbound}
\nonumber
&\max_{1\leq r \leq k} \pr(Z_{m}=r | Z_{0}=j)\\
&\hspace{1cm}\leq e^{Q_{m}r}
\sum_{l=1}^{j} \binom{j}{l}Q_{m}^{l}(1-Q_{m})^{j-l}e^{-l/2}\max_{r\geq 1}\widehat{\pr}(S_{l}^{(m)}=r).  
\end{align}
Combining \eqref{eq:exponentialchange} and \eqref{eq:inequalityforexponetial}, we also have 
$$
\sup_{r\geq 1}\widehat{\pr}(\xi_{1}^{(m)} =r) 
\leq e\sup_{r \geq 1}\pr(\xi_{1}^{(m)}=r) \leq \frac{C}{m},
$$
where in the last step we have used Proposition from \cite{NV06}. Using now the fact that the concentration function of the sum is smaller than the concentration function of a single summand, we have
$$
\sup_{r\geq 1}\widehat{\pr}(S_{l}^{(m)}=r )\leq \frac{C}{m}.
$$
Plugging this into \eqref{mainbound} and using the elementary bound $1-x\le e^{-x}$, we obtain
\begin{align*}
\max_{1\leq r \leq k}\pr(Z_{m}=r | Z_{0}=j) 
&\leq C\frac{e^{Q_{m}k}}{m}\sum_{l=1}^{j}\binom{j}{l}Q_{m}^{l}(1-Q_{m})^{j-l}e^{-l/2}\\
& =C\frac{e^{Q_{m}k}}{m} (1-Q_{m}+Q_{m}e^{-1/2})^{j}\\
&\le \frac{C}{m}e^{Q_{m}(k-(1-e^{-1/2})j)}.
\end{align*}

Recalling that $Q_{m}\sim \frac{2}{Bm}$ as $m\to \infty$, we conclude that
\begin{equation}\label{eq:max-bound}
\max_{1\leq r \leq k}\pr(Z_{m}=r | Z_{0}=j)\leq CQ_{m}e^{-\frac{(1-e^{-1/2})}{2}jQ_{m}},    
\end{equation}
for all $m\geq m_{0}$ and all $j\geq \frac{ 2k}{(1-e^{-1/2})}$.
The assumption $m=o(k)$ implies that $kQ_m\to\infty$ and, consequently,
$$
\max_{1\leq r \leq k,j\ge ck}\pr(Z_{m}=k | Z_{0}=j)
=o\left(\frac{1}{k} \right)
$$
for every $c>2(1-e^{-1/2})^{-1}$. The second claim is also immediate from \eqref{eq:max-bound}.
\end{proof}
\begin{lemma}
There exists a constant $C$ such that 
\begin{align}
\label{eq:app.7}
\sup_{l\ge 1}\pr(Z_m=l|Z_0=j)\le\frac{C}{\sqrt{mj}}
\end{align}   
for all $m,j\ge1$.
\end{lemma}
\begin{proof}
 The following inequality has been derived in the proof of Lemma 20 in \cite{NV06}:
$$
\sup_{l\ge 1}\pr(Z_m=l|Z_0=j)
\le\frac{C}{m}\e\left[\zeta^{-1/2}_{j,m};\zeta_{j,m}>0\right],
$$
where $\zeta_{j,m}$ has binomial distribution with parameters $j$ and $\pr(Z_m>0)$. Combining this with Lemma 18 from \cite{NV06}, we obtain the desired bound. 
\end{proof}
\subsection{Local estimates for the process with immigration.}
It has been shown in \cite{ShW} that there exists a constant $C$ such that 
\begin{align}
\label{eq:Ubound-1}
\pr(Y_n=l|Y_0=0)\le \frac{C}{l},\quad l\ge1.
\end{align}
Comparing this with the local limit theorem, we infer that this upper bound has the right order for $l$ comparable with $n$. The main purpose of this subsection is to derive a bit better bound in the region $l\gg n$ of large deviations. More precisely, we want to show that if $l/n\to\infty$ then
$$
l\,\pr(Y_n=l|Y_0=0)\to0.
$$

To keep the minimal moment restrictions we shall use truncation of random variables describing the immigration component of the process $\{Y_n\}$. 

Let $\{\widetilde{\eta}_n\}$ be a sequence of i.i.d. random variables with the following distribution 
$$
\pr\big(\widetilde\eta_1=l\big)
:=
\pr\big(\eta_1=l \mid \eta_1\le N\big),\quad  l\ge0,
$$
where $N$ is an arbitrary natural number.
Let $\{\widetilde{Y}_n\}$ be a critical branching process with offspring random variables $\{\xi_{n,j}\}$ and with immigration variables $\{\widetilde\eta_n\}$, that is,
$$
\widetilde{Y}_{n+1}=\sum_{j=1}^{\widetilde{Y}_{n}} \xi_{n+1,j}+\widetilde{\eta}_{n+1}, \quad n \geq 0.
$$
We shall always assume that $\widetilde{Y}_0=0$.

Our next goal will be to derive an appropriate upper bound for the probability $\pr(\widetilde{Y}_n=l)$.
We start by comparing the generating functions of $\eta_1$ and $\widetilde\eta_1$.
Define
$$
\widetilde h(s):=\mathbf{E} s^{\widetilde\eta_1}, \ \ |s|\leq 1
$$
and
$$
h(s):=\mathbf{E} s^{\eta_1}=\sum_{k=0}^\infty q_ks^k, \ \ |s|\leq 1.
$$
\begin{lemma}
\label{lem:h-tranc}    
Assume that $N$ is such that $\pr(\eta_1>N)\le 1/4$. Then, for all
$s$ such that $|s|\le1$ and $|h(s)|\ge 1/2$,
$$
\left|\log\widetilde h(s)-\log h(s)\right|\le 5\pr(\eta_1>N).
$$
$\log$ denotes here the principle value of the logarithm.
\end{lemma}
\begin{proof}
We first notice that the equality
$$
\overline{h}(s):=\sum_{l=0}^Nq_ls^l=\pr(\eta_1\le N)\widetilde h(s)
$$
implies that 
$$
\left|\log \overline h(s)-\log\widetilde h(s)\right|
\le -\log\left(1-\pr(\eta_1>N)\right).
$$
Decomposing the right hand side into the Taylor series and using the assumption $\pr(\eta_1>N)\le 1/4$, we conclude that
\begin{align}
\label{eq:trunc-h.1}  
\left|\log \overline h(s)-\log\widetilde h(s)\right|
\le \frac{5}{4}\pr(\eta_1>N).
\end{align}
Next, 
\begin{align*}
\log\overline h(s)
&=\log \left(h(s)-(h(s)-\overline h(s))\right)\\
&=\log h(s)+\log\left(1-\frac{h(s)-\overline h(s)}{h(s)}\right).
\end{align*}
The assumptions $|h(s)|\ge 1/2$ and $\pr(\eta_1>N)\le 1/4$ imply that
$$
\left|\frac{h(s)-\overline h(s)}{h(s)}\right|
\le \frac{\pr(\eta_1>N)}{|h(s)|}
\le 2\pr(\eta_1>N)
\le \frac{1}{2}.
$$
Using once again the Taylor series for $\log(1-x)$, we obtain 
$$
\left|\log\overline h(s)-\log h(s)\right|
\le \frac{3}{2}\left|\frac{h(s)-\overline h(s)}{h(s)}\right|
\le 3\pr(\eta_1>N).
$$
Combining this with \eqref{eq:trunc-h.1}, we obtain the desired bound.
\end{proof}
Set
$$
\widetilde{H}_{n}(s)
:=
\mathbf E\!\left[s^{\widetilde Y_{n}}\right]
=
\prod_{j=0}^{n-1}\widetilde h\!\bigl(f_{j}(s)\bigr),\quad |s|\le 1,
$$
where $f_j(s)$ denotes the $j$-th iteration of the generating function $f(s)$.
\begin{lemma}
\label{lem:trunc-H}
Assume that $\pr(\eta_1>N)\le 1/4$. Then there exists an independent of $N$ constant $C$ such that, for all 
$n\ge 1$ and $|t|\in[1/n,\pi/2]$,
$$
|\widetilde H_n(e^{it})|
\le\left|\frac{\widetilde H_n(e^{it})}{\widetilde h(e^{it})}\right|
\le C e^{5n\pr(\eta_1>N)} (n|t|)^{-\theta},
$$
where $\theta=\min\{\gamma,1/2\}$.
\end{lemma}
\begin{proof}
Applying (35) from \cite{ShW} to the function  $\widetilde H_n$, we obtain 
$$
\left|\frac{\widetilde H_n(s)}{\widetilde h(s)}\right|
\le \exp\left\{\sum_{j=N}^{n-1}{\rm Re}\log\widetilde h(f_j(s))\right\},
$$
where ${\rm Re}\,z$ denotes the real part of $z$.

Since $\lim_{j\to\infty}f_j(s)=1$ uniformly in $s$ from the unit disc, there exists $M$ such that $|h(f_j(s))|\ge 1/2$ for all $j\ge M$. This allows us to apply Lemma~\ref{lem:h-tranc}, which gives 
$$
\left|\frac{\widetilde H_n(s)}{\widetilde h(s)}\right|
\le \exp\left\{5n\pr(\eta_1>N)+
\sum_{j=M}^{n-1}{\rm Re}\log
h(f_j(s))\right\}.
$$
Having this inequality, we may apply Lemmata 7 and 8 from \cite{ShW}
and obtain the desired bound.
\end{proof}
\begin{lemma}\label{lem:trunc-Y}
Assume that $N$ is such that $n\pr(\eta_1>N)\le 1$. Then
there exists a constant $C$ such that, for all $n,l\geq 1$,
$$
\pr(\widetilde{Y}_{n}=l)\leq \frac{Cn}{l^{2}}+\frac{C\widetilde{h}''(1)}{l^{2}}.  
$$
\end{lemma}
\begin{proof}
 For $l = 1$ the bound is trivial, so we assume $l\ge 2$. We use the following bound for the concentration function (see \cite{Petrov75}):
 \begin{equation}\label{eq:concformula}
     \sup_{x}\pr(\xi=x)\leq \left(\frac{96}{95} \right)^{2}\frac{1}{a}\int_{-a}^{a}\left|\phi(t)\right|\mathrm{d}t,
 \end{equation}
where $\phi(t)$ is the characteristic function of $\xi$ and $a$ is an arbitrary positive number.
Let $\xi$ be a random variable with the generating function $\widetilde{H}_{n}''(s)/\widetilde{H}_{n}''(1)$. Then there exists a constant $c>0$ such that, for any $n\geq 1$,
\begin{equation}\label{eq:concformula1}
     \sup_{l\geq 2}l(l-1)\pr(\widetilde{Y}_{n}=l)\leq \frac{c}{a}\int_{-a}^{a}\left|\widetilde{H}_{n}''(e^{\mathrm{i}t})\right|\mathrm{d}t.
 \end{equation}
Thus, it suffices to show that
\begin{equation}\label{eq:mainbound}
\int_{-a}^{a} \left| \widetilde{H}_{n}''\left(e^{\mathrm{i}t}\right) \right|\,\mathrm{d}t \leq Cn+C\widetilde{h}''(1).
\end{equation}
For this purpose, we first observe that, for $|s|\leq 1$,
$$
\widetilde{H}_{n}'(s)
=
\widetilde{H}_{n}(s)\sum_{l=0}^{n-1}
\frac{\widetilde{h}'\bigl(f_{l}(s)\bigr)}{\widetilde{h}\bigl(f_{l}(s)\bigr)}
f_{l}'(s).
$$
Differentiating this equality, we obtain
$$
\widetilde{H}_{n}''(s)=\widetilde{H}_{n}'(s)\sum_{l=0}^{n-1}
\frac{\widetilde{h}'\bigl(f_{l}(s)\bigr)}{\widetilde{h}\bigl(f_{l}(s)\bigr)}f_{l}'(s)+\widetilde{H}_{n}(s)\sum_{l=0}^{n-1}\frac{\mathrm{d}}{\mathrm{d}s}\left(
\frac{\widetilde{h}'\bigl(f_{l}(s)\bigr)}{\widetilde{h}\bigl(f_{l}(s)\bigr)}f_{l}'(s) \right)
$$
$$
=\widetilde{H}_{n}(s)\left(\sum_{l=0}^{n-1}
\frac{\widetilde{h}'\bigl(f_{l}(s)\bigr)}{\widetilde{h}\bigl(f_{l}(s)\bigr)}f_{l}'(s)  \right)^{2}
$$
$$
+\widetilde{H}_{n}(s)\sum_{l=0}^{n-1}\left(\frac{\widetilde{h}'\bigl(f_l(s)\bigr)}{\widetilde{h}\bigl(f_{l}(s)\bigr)}f_{l}''(s)+\frac{\widetilde{h}''\bigl(f_{l}(s)\bigr)}{\widetilde{h}\bigl(f_{l}(s)\bigr)}(f_{l}'(s))^{2}-\left(\frac{\widetilde{h}'\bigl(f_{l}(s)\bigr)}{\widetilde{h}\bigl(f_{l}(s)\bigr)}\right)^{2}(f_{l}'(s))^{2} \right)
$$
$$
=\widetilde{H}_{n}(s)\sum_{l=0}^{n-1}\left(\frac{\widetilde{h}'\bigl(f_{l}(s)\bigr)}{\widetilde{h}\bigl(f_{l}(s)\bigr)}f_{l}''(s)+\frac{\widetilde{h}''\bigl(f_{l}(s)\bigr)}{\widetilde{h}\bigl(f_{l}(s)\bigr)}(f_{l}'(s))^{2} \right)
$$
$$
+2\widetilde{H}_{n}(s)\sum_{l=0}^{n-1}\sum_{j=0}^{l-1}\frac{\widetilde{h}'\bigl(f_{l}(s)\bigr)}{\widetilde{h}\bigl(f_{l}(s)\bigr)}f_{l}'(s)\frac{\widetilde{h}'\bigl(f_{j}(s)\bigr)}{\widetilde{h}\bigl(f_{j}(s)\bigr)}f_{j}'(s).
$$

Notice next that 
\begin{align*}
|\widetilde{H}_{n}''(e^{\mathrm{i}t})| 
\leq \widetilde{H}_{n}''(1)
&=\sum_{l=0}^{n-1}\left(\widetilde h'(1)f_l''(1)+\widetilde h''(1)\right)
+2\sum_{l=0}^{n-1}\sum_{j=0}^{l-1}\left(\widetilde h'(1)\right)^2\\
&\le Cn(n+\widetilde h''(1)).
\end{align*}
This implies that
\begin{align}
\label{eq:int-smaller}
\int_{|t|<1/n}|\widetilde{H}_{n}''(e^{\mathrm{i}t})|\mathrm{d}t \leq \frac{2}{n}\widetilde{H}_{n}''(1)\leq C(n+\widetilde h''(1)).
\end{align}

To estimate the integral of $\widetilde{H}_{n}''(e^{\mathrm{i}t})$
over $\{|t|\in(1/n,a)\}$ we list necessary estimates for the summands in the representation for $\widetilde{H}_{n}''(s)$. It is easily verified that 
$$
f_{n}''(s)=f''(f_{n-1}(s))(f_{n-1}'(s))^{2}+f'(f_{n-1}(s))f_{n-1}''(s).
$$
Recalling that $|f''(s)|\le B$ and $|f'(s)|\le 1$ for all $|s|\le 1$, we conclude that 
$$
|f_{n}''(s)|\le B|f_{n-1}'(s)|^{2}+|f_{n-1}''(s)|.
$$
Iterating this bound, we infer that for every $l\in\{1,2,\ldots,n-1\}$ one has the inequality
\begin{equation}
\label{eq:2.der-prelim}
|f_n''(s)|\le |f_l''(s)|+B\sum_{j=l+1}^{n-1}|f_{j}'(s)|^{2}
\end{equation}
Set $s=e^{\mathrm{i}t}$. It is known (see page 409 in \cite{NV06}) that 
\begin{equation}\label{eq:boundforder}
\left|f_{j}'(e^{\mathrm{i}t})\right| \leq \frac{C}{(j|t|)^{3/2}},\quad |t|\geq \frac{1}{j}.
\end{equation}
Notice also that
\begin{equation}\label{eq:boundforsecder}
|f_{l}''(e^{\mathrm{i}t})| \leq Bl, \ \ t\in \mathbb{R}.
\end{equation}
Applying \eqref{eq:boundforder} and \eqref{eq:boundforsecder} to the right hand side of \eqref{eq:2.der-prelim}, we conclude that
$$
|f_n''(e^{\mathrm{i}t})|\le Bl+C\sum_{j=l+1}^{n-1}\frac{1}{|t|^{3}j^{3}}, 
\quad |t|\ge\frac{1}{l+1}.
$$
Choosing $l$ so that $1/(l+1)\le |t|< 1/l$, we obtain
\begin{align}\label{eq:boundforsecder1}
\nonumber
|f_{n}''(e^{\mathrm{i}t})| 
&\leq C_1B\left(l+\sum_{j=l+1}^{n-1}\frac{1}{|t|^{3}j^{3}}\right)\\
&\leq C_1B\left(l+\frac{1}{l^{2}|t|^{3}} \right)
\leq C_2l \leq \frac{C_3}{|t|}.
\end{align}

It is clear that $|\widetilde{h}'(f_{l}(s))|\leq \widetilde{h}'(1)$ for all $s$ in the unit disc. It is well-known that $f_{l}(s)\to 1$ as $l\to\infty$ uniformly in $|s|\leq 1$. This property implies that, for sufficiently small $a>0$, 
$$ \left|\widetilde h\left(f_{l}\left(e^{\mathrm{i}t}\right)\right)\right| \geq \frac{1}{2} $$
for all $l\geq 1$ and all $t\in[-a,a]$. 
Combining the statements above, we note that $\left|\frac{\widetilde{h}'(f_{l}(e^{\mathrm{i}t}))}{\widetilde{h}(f_{l}(e^{\mathrm{i}t}))} \right|$ is bounded.

Combining Lemma~\ref{lem:trunc-H} and \eqref{eq:boundforsecder}, we conclude that
$$
\begin{aligned}
\int_{1/n<|t|<1/l}|\widetilde{H}_{n}(e^{\mathrm{i}t})f_{l}''(e^{\mathrm{i}t})|\mathrm{d}t 
& \leq 
C_1l \int_{1/n<|t|<1/l}|\widetilde{H}_{n}(e^{it})|\mathrm{d}t \\
& \leq C_2l \int_{1/n<|t|<1/l} \frac{\mathrm{d}t}{(n|t|)^{\theta}} \le \frac{C_3l^{\theta}}{n^{\theta}}.
\end{aligned}
$$
If $1/l<|t|<a$, then we combine \eqref{eq:boundforsecder1} and Lemma~\ref{lem:trunc-H} to obtain
$$
\int_{1/l<|t|<a}|\widetilde{H}_{n}(e^{\mathrm{i}t})f_{l}''(e^{\mathrm{i}t})|\mathrm{d}t \leq 
\frac{C_1}{n^{\theta}} \int_{1/l<|t|<a}\frac{1}{|t|^{\theta+1}}\mathrm{d}t \le \frac{C_2l^{\theta}}{n^{\theta}}.
$$
Therefore,
\begin{align}
\label{eq:int-bigger}
\sum_{l=0}^{n-1}
\int_{1/n<|t|<a}|\widetilde{H}_{n}(e^{\mathrm{i}t})f_{l}''(e^{\mathrm{i}t})|\mathrm{d}t 
\le C_1\sum_{l=1}^{n-1}\frac{l^\theta}{n^\theta}
\le C_2 n.
\end{align}
Recalling that the functions $\left|\frac{\widetilde{h}'(f_{l}(e^{\mathrm{i}t}))}{\widetilde{h}(f_{l}(e^{\mathrm{i}t}))} \right|$
are bounded, we infer from \eqref{eq:int-bigger} that 
\begin{align}
\label{eq:int-bigger.1}
\sum_{l=0}^{n-1}
\int_{1/n<|t|<a}\left|\widetilde{H}_{n}(e^{\mathrm{i}t})f_{l}''(e^{\mathrm{i}t})\frac{\widetilde{h}'(f_{l}(e^{\mathrm{i}t}))}{\widetilde{h}(f_{l}(e^{\mathrm{i}t}))}\right|\mathrm{d}t 
\le C n.
\end{align}

Next, using first the bounds $|f_{l}'(e^{\mathrm{i}t})|\leq 1$, $|\widetilde{h}''(f_{l}(e^{\mathrm{i}t}))|\leq |\widetilde{h}''(1)|$ and applying then Lemma~\ref{lem:trunc-H}, we obtain
\begin{align*}
&\int_{1/n<|t|<a}\left|
\widetilde{H}_{n}(e^{\mathrm{i}t})\frac{\widetilde{h}''(f_{l}(e^{\mathrm{i}t}))}{\widetilde{h}(f_{l}(e^{\mathrm{i}t}))} (f_{l}'(e^{\mathrm{i}t}))^{2}\right|\mathrm{d}t \\
&\hspace{1cm} \leq \widetilde{h}''(1)\int_{1/n<|t|<a}|\widetilde{H}_{n}(e^{\mathrm{i}t})f_{l}'(e^{\mathrm{i}t})|\mathrm{d}t \\
&\hspace{1cm} \leq C \widetilde{h}''(1)
\left(\int_{1/n<|t|<1/l}(n|t|)^{-\theta}\mathrm{d}t
+\int_{1/l<|t|<a}(n|t|)^{-\theta}|f_{l}'(e^{\mathrm{i}t})|\mathrm{d}t \right).
\end{align*}
Applying now \eqref{eq:boundforder}, we conclude that 
$$
\int_{1/n<|t|<a}\left|
\widetilde{H}_{n}(e^{\mathrm{i}t})\frac{\widetilde{h}''(f_{l}(e^{\mathrm{i}t}))}{\widetilde{h}(f_{l}(e^{\mathrm{i}t}))} (f_{l}'(e^{\mathrm{i}t}))^{2}\right|\mathrm{d}t
\le C\widetilde h''(1)\frac{l^{\theta-1}}{n^\theta}.
$$
Consequently,
\begin{align}
\label{eq:int-bigger.2}
\sum_{l=0}^{n-1}\int_{1/n<|t|<a}\left|
\widetilde{H}_{n}(e^{\mathrm{i}t})\frac{\widetilde{h}''(f_{l}(e^{\mathrm{i}t}))}{\widetilde{h}(f_{l}(e^{\mathrm{i}t}))} (f_{l}'(e^{\mathrm{i}t}))^{2}\right|\mathrm{d}t \leq C\widetilde{h}''(1).
\end{align}
Finally, we split the integral
$$
\int_{1/n < |t|<a}|\widetilde{H}_{n}(e^{\mathrm{i}t})f_{i}'(e^{\mathrm{i}t})f_{j}'(e^{\mathrm{i}t})|\mathrm{d}t
=I_{1}+I_{2}+I_{3},
$$
where
\begin{align*}
I_{1}&:=\int_{1/n < |t|<1/l}|\widetilde{H}_{n}(e^{\mathrm{i}t})f_{i}'(e^{\mathrm{i}t})f_{j}'(e^{\mathrm{i}t})|\mathrm{d}t,\\
I_{2}&:=\int_{1/l < |t|<1/j}|\widetilde{H}_{n}(e^{\mathrm{i}t})f_{i}'(e^{\mathrm{i}t})f_{j}'(e^{\mathrm{i}t})|\mathrm{d}t,\\
I_{3}&:=\int_{1/j < |t|<a}|\widetilde{H}_{n}(e^{\mathrm{i}t})f_{i}'(e^{\mathrm{i}t})f_{j}'(e^{\mathrm{i}t})|\mathrm{d}t.
\end{align*}

We now estimate $I_1$, $I_2$ and $I_3$ separately.

Using the trivial bounds $|f_{l}'(e^{\mathrm{i}t})|\leq 1$,
$|f_{j}'(e^{\mathrm{i}t})|\leq 1$ and Lemma~\ref{lem:trunc-H}, it yields
$$
I_{1}\leq \int_{1/n < |t|<1/l}|\widetilde{H}_{n}(e^{\mathrm{i}t})|\mathrm{d}t \leq C_1\int_{1/n < |t|<1/l}\frac{1}{(n|t|)^{\theta} }\mathrm{d}t\leq \frac{C_2l^{\theta-1}}{n^{\theta}}.
$$
To estimate $I_2$, we replace $|f_{j}'(e^{\mathrm{i}t})|$ by $1$ and apply then Lemma~\ref{lem:trunc-H} and \eqref{eq:boundforder}:
$$
I_{2} 
\leq \int_{1/l < |t|<1/j}|\widetilde{H}_{n}(e^{\mathrm{i}t})
f_l(e^{\mathrm{i}t})|\mathrm{d}t
\leq C_1\int_{1/l < |t|<1/j}\frac{1}{(n|t|)^{\theta}}\frac{1}{(lt)^{3/2}}\mathrm{d}t \le \frac{C_2l^{\theta-1}}{n^{\theta}}.
$$
Applying \eqref{eq:boundforder} to $f'_l$ and to $f'_j$, and using for the third time Lemma~\ref{lem:trunc-H}, we conclude that
$$
\begin{aligned}
I_{3} & \leq C_1\int_{1/j < |t|<a}\frac{1}{(n|t|)^{\theta}} \frac{1}{(l|t|)^{3/2}} \frac{1}{(j|t|)^{3/2}}\mathrm{d}t \\
& =\frac{C_1}{n^{\theta}l^{3/2}j^{3/2}}\int_{1/j<|t|<a}\frac{1}{|t|^{3+\theta}}\mathrm{d}t\leq \frac{C_2j^{\theta+1/2}}{n^{\theta}l^{3/2}}.
\end{aligned}
$$
Gathering these bounds, we obtain
$$
\int_{1/n < |t|<a}|\widetilde{H}_{n}(e^{\mathrm{i}t})f_{i}'(e^{\mathrm{i}t})f_{j}'(e^{\mathrm{i}t})|\mathrm{d}t
\le C\left(\frac{l^{\theta-1}}{n^\theta}
+\frac{j^{\theta+1/2}}{n^{\theta}l^{3/2}}\right)
$$
and, consequently,
\begin{align}
\label{eq:int-bigger.3}
\nonumber
&\sum_{l=0}^{n-1}\sum_{j=1}^{l-1}\int_{1/n < |t|<a}|\widetilde{H}_{n}(e^{\mathrm{i}t})f_{i}'(e^{\mathrm{i}t})f_{j}'(e^{\mathrm{i}t})|
\mathrm{d}t\\
&\hspace{2cm}
\le C_1\sum_{l=0}^{n-1}\sum_{j=1}^{l-1}\left(\frac{l^{\theta-1}}{n^\theta}
+\frac{j^{\theta+1/2}}{n^{\theta}l^{3/2}}\right)
 \leq {C_2n}.
\end{align}
From \eqref{eq:int-bigger.1}, \eqref{eq:int-bigger.2} and \eqref{eq:int-bigger.3} we infer that 
$$
\int_{1/n<|t|<a}|\widetilde{H}_{n}''(e^{\mathrm{i}t})|\mathrm{d}t 
\leq
C(n+\widetilde h''(1)).
$$
Combining this with \eqref{eq:int-smaller}, we finish the proof.
\end{proof}
\begin{corollary}
If $\e\eta_1^2<\infty$ then there exists a constant $C$ such that
$$
\pr(Y_n=l)\le C\frac{n}{l^2},
\quad l,n\ge1.
$$
\end{corollary}
\begin{proof}
 Set
 $$
 A_{N,n}:=\left\{\eta_k\le N\text{ for all }k\le n\right\}.
 $$
 Then, for all fixed $l,n\ge 1$,
 \begin{align*}
  \pr(Y_n=l)
  &=\lim_{N\to\infty}\pr(Y_n=l,A_{N,n})
  =\lim_{N\to\infty}\pr(\widetilde Y_n=l)\pr(A_{N,n})\\
  &=\lim_{N\to\infty}\pr(\widetilde Y_n=l).
 \end{align*}
 Applying now Lemma~\ref{lem:trunc-Y}, we conclude that 
 $$
 \pr(Y_n=l)\le C\lim_{N\to\infty}\frac{n+\widetilde{h}''(1)}{l^2}
 \le C \frac{n+\e\eta_1^2}{l^2}.
 $$
 Thus, the proof is complete.
\end{proof}
If one does not assume the existence of the second moment of $\eta_1$ then one has the following bound.
\begin{lemma}
\label{lem:loc-ld}
There exist a function $\delta(l)=o(l)$ and a constant $C$ such that 
$$
\pr(Y_n=l)\le C\frac{n+\delta(l)}{l^2}
\quad\text{for all }l\ge 2n\ge1.
$$
\end{lemma}
\begin{proof}
We have the following decomposition:
\begin{align*}
\pr(Y_n=l)
&=\pr(Y_n=l,A_{N,n})+\pr(Y_n=l,A^c_{N,n})\\
&\le \pr(\widetilde Y_n=l)+\pr(Y_n=l,A^c_{N,n}).
\end{align*}
We choose here $N=al$, where $a$ is a constant which will be chosen later. The assumption $\e\eta_1<\infty$ implies that
$\widetilde h''(1)=o(l)$ as $l\to\infty.$ Therefore,
by Lemma~\ref{lem:trunc-Y},
$$
\pr(Y_n=l)
\le C\frac{n+\delta_1(l)}{l^2}+\pr(Y_n=l,A^c_{N,n}),
$$
where $\delta_1(l)=o(l)$.

For the probability term on the right hand side we have 
$$
\pr(Y_n=l,A^c_{N,n})
\le \sum_{k=1}^{n-1}\pr(\eta_k>N)\max_{j>N}\pr(Y_{n-k}=l|Y_0=j).
$$
Furthermore, by the total probability law,
\begin{align*}
\pr(Y_{n-k}=l|Y_0=j)
&=\sum_{r=0}^l\pr(Y_{n-k}=r|Y_0=0)\pr(Z_{n-k}=l-r|Z_0=j)\\
&\le \max_{r\ge l/2}\pr(Y_{n-k}=r|Y_0=0)
+\max_{r\le l/2}\pr(Z_{n-k}=l-r|Z_0=j).
\end{align*}
According to \eqref{eq:Ubound-1},
$$
\max_{r\ge l/2}\pr(Y_{n-k}=r|Y_0=0)\le \frac{C}{l}.
$$
Moreover, if $a>\frac{2}{(1-e^{-1/2})}$ then,
by Lemma~\ref{lem:max-bound},
$$
\max_{r\le l/2}\pr(Z_{n-k}=l-r|Z_0=j)\le \frac{C}{l}.
$$
Therefore,
$$
\pr(Y_n=l,A^c_{N,n})\le \frac{C}{l}n\pr(\eta_1>N)
\le C\frac{n}{l^2}.
$$
This finishes the proof of the lemma.
\end{proof}
\section{Approximation by a process with a geometric offspring law}
For every $m\ge1$ we define
$$
a_m(s)=1-\frac{1}{\frac{1}{1-s}+\frac{Bm}{2}}.
$$
It is easy to see that $a_m(s)=a_{m-1}(a_1(s))$ for every $m\ge2$ and that $a_1'(1)=1$. Therefore, this sequence of generating function corresponds to a critical Galton--Watson process with the offspring generating function $a_1(s)$. We shall denote this process by $\{G_m\}$. 
The main purpose of this section is to compare mass functions of random variables $Z_m$ and $G_m$ conditioned on $Z_0=G_0=j\ge1$.
\begin{proposition}
\label{prop:comparison}
There exists a function $\beta(x)\to0$ as $x\to\infty$ such that, for all $j\ge 1$,
$$
\sup_{l\ge m}\Big|\pr(Z_m=l|Z_0=j)-\pr(G_m=l|G_0=j)\Big|
\le \frac{\beta(m)}{m^{3/2}}\sqrt{j}.
$$
\end{proposition}
\begin{proof}
According to formula (2.7) in \cite{NV06},
\begin{align}
\label{eq:app.0}
f_m(s)-a_m(s)=d_m(s)g_m(s), 
\end{align}
where
$$
d_m(s):=\frac{1}{1-f_m(s)}-\frac{1}{1-a_m(s)}
=\frac{1}{1-f_m(s)}-\frac{1}{1-s}-\frac{Bm}{2}
$$
and 
$$
g_m(s):=(1-f_m(s))(1-a_m(s))
=(1-f_m(s))\left(\frac{1}{1-s}+\frac{Bm}{2}\right).
$$
By Lemmata 2 and 21 from \cite{NV06},
\begin{align}
\label{eq:app.1}
\sum_{k=0}^\infty |[s^k]d_m(s)|=o(m)
\quad\text{as } m\to\infty
\end{align}
and 
\begin{align}
\label{eq:app.2}
\lim_{k\to\infty}[s^k]d_m(s)=0
\quad\text{uniformly in }m\ge1.
\end{align}
Here and in what follows, $[s^k]r(s)$ denotes the $k$-th coefficient of a generating function $r(s)$.

Furthermore, by Lemma 19 from the same paper,
\begin{align}
\label{eq:app.3}
[s^0]g_m(s)\le\frac{C}{m^2}
\quad\text{and}\quad 
\sup_{k\ge1}\Big|[s^k]g_m(s)\Big|\le\frac{C}{m^3}
\end{align}
with some independent of $m$ constant $C<\infty$.

Combining \eqref{eq:app.1} -- \eqref{eq:app.3}, we conclude that 
\begin{align}
\label{eq:app.4}
\nonumber
\Big|[s^l]d_m(s)g_m(s)\Big|
&\le \Big|[s^l]d_m(s)\Big|g_m(0)
+\sum_{k=1}^l \Big|[s^k]g_m(s)\Big|\Big|[s^{l-k}]d_m(s)\Big|\\
\nonumber 
&\le \frac{C}{m^2}\Big|[s^l]d_m(s)\Big|
+\frac{C}{m^3}\sum_{j=0}^{l-1} \Big|[s^j]d_m(s)\Big|\\
&\le \frac{\beta_1(l)+\beta_2(m)}{m^2},
\end{align}
where $\beta_1(x),\beta_2(x)\to0$ as $x\to\infty$. 

To estimate the difference $\pr(Z_m=l|Z_0=j)-\pr(G_m=l|G_0=j)$, we write down the representation
\begin{align*}
&\pr(Z_m=l|Z_0=j)-\pr(G_m=l|G_0=j) \\
&\hspace{1cm} =[s^l](f_m^j(s)-a_m^j(s))
=[s^l]\left((f_m(s)-a_m(s))\sum_{t=0}^{j-1}f_m^t(s)a_m^{j-1-t}(s)\right)\\
&\hspace{1cm}=\sum_{k=0}^l
[s^k]h_m(s)g_m(s)
[s^{l-k}]\left(\sum_{t=0}^{j-1}f_m^t(s)a_m^{j-1-t}(s)\right).
\end{align*}
Using \eqref{eq:app.4} and the relations
\begin{align}
\label{eq:app.4.5}
1-f_m(0)\sim\frac{2}{Bm},\quad  1-a_m(0)\sim\frac{2}{Bm}
\quad\text{as }m\to\infty,
\end{align}
we obtain 
\begin{align}
\label{eq:app.6}
\nonumber 
\left|[s^l]d_m(s)g_m(s)
[s^{0}]\left(\sum_{t=0}^{j-1}f_m^t(s)a_m^{j-1-t}(s)\right)\right|
&\le \frac{\beta_1(l)+\beta_2(m)}{m^2}
\sum_{t=0}^{j-1}f_m^t(0)a_m^{j-1-t}(0)\\
&\le
\frac{\beta_1(l)+\beta_2(m)}{m^2} je^{-\theta j/m}
\end{align}
with some $\theta>0$.

Applying \eqref{eq:app.7} to the processes $\{Z_m\}$ and $\{G_m\}$, we conclude that, uniformly in $t\in\{0,1,\ldots,j-1\}$,
$$
\sup_{l\ge 1}[s^l]\left(f_m^t(s)a_m^{j-1-t}(s)\right)
\le \frac{C}{\sqrt{mj}}.
$$
Consequently,
\begin{align*}
&\left|\sum_{k=0}^{l-1}
[s^k]d_m(s)g_m(s)
[s^{l-k}]\left(\sum_{t=0}^{j-1}f_m^t(s)a_m^{j-1-t}(s)\right)\right|\\ &\hspace{1cm}\le \frac{C}{\sqrt{mj}}j\sum_{k=0}^{l-1}\Big|[s^k]d_m(s)g_m(s)\Big|\\
&\hspace{2cm}\le \frac{C}{\sqrt{mj}}j
\left(\sum_{k=0}^\infty\Big|[s^k]d_m(s)\Big|\right)
\left(\sum_{k=0}^\infty\Big|[s^k]g_m(s)\Big|\right).
\end{align*}
It is immediate from the definition of $g_m$ and from \eqref{eq:app.4.5} that 
$$
\sum_{k=0}^\infty\Big|[s^k]g_m(s)\Big|\le\frac{C}{m^2}.
$$
From this and from \eqref{eq:app.1} we infer that 
\begin{align}
\label{eq:app.8}
\left|\sum_{k=0}^{l-1}
[s^k]d_m(s)g_m(s)
[s^{l-k}]\left(\sum_{t=0}^{j-1}f_m^t(s)a_m^{j-1-t}(s)\right)\right|
\le \frac{\beta_3(m)}{m^{3/2}}\sqrt{j},
\end{align}
where $\beta_3(x)\to0$ as $x\to\infty$. Combining \eqref{eq:app.6} and \eqref{eq:app.8}, we conclude that
\begin{align*}
&\sup_{l\ge m}\Big|\pr(Z_m=l|Z_0=j)-\pr(G_m=l|G_0=j)\Big|\\
&\hspace{1cm}
\le \frac{\sup_{l\ge m}\beta_1(l)+\beta_2(m)}{m^2} je^{-\theta j/m}
+\frac{\beta_3(m)}{m^{3/2}}\sqrt{j}\\
&\hspace{2cm}\le\left(C_\theta\left(\sup_{l\ge m}\beta_1(l)+\beta_2(m)\right)+\beta_3(m)\right)\frac{\sqrt{j}}{m^{3/2}}.
\end{align*}
This completes the proof of the proposition.
\end{proof}
\section{Local probabilities of a process with a geometric offspring law.}

We first derive an upper bound for the local probability
$\pr(G_m=l|G_0=j)$.
\begin{proposition}\label{prop:geom-ub}
There exists a constant $C$ such that
$$
\pr(G_m=l|G_0=j)\le \frac{Cj^{1/4}}{m^{1/2}l^{3/4}}
\exp\left\{-\frac{(l-j)^2}{2Bml}\right\}
$$
for all $j\le l\le m^2$.
\end{proposition}
\begin{proof}
As it is usual for large deviation probabilities, we shall perform the exponential change of measure. Let $s_0>1$ be such that $a_m(s_0)<\infty$. For every such $s_0$ we may define
$$
\widehat{a}_m(s)=\frac{a_m(s_0s)}{a_m(s_0)}.
$$
The optimal choice of the parameter $s_0$ is the solution to the equation
$$
\widehat{a}_m'(1)=\frac{s_0a_m'(s_0)}{a_m(s_0)}=\frac{l}{j}.
$$
For our purposes it will be sufficient to take an 'approximative' solution
$$
s_0=1+\frac{2}{Bm}\left(1-\sqrt{\frac{j}{l}}\right)>1.
$$
For each $m\ge1$, let $\{U^{(m)}_j\}_{j\ge 1}$ be independent, identically distributed random variables with the generating function $\widehat{a}_m(s)$. Then one has 
\begin{align}
\label{eq:B.1}
\pr(G_m=l|G_0=j)
=s_0^{-l}a_m^j(s_0)\pr(U^{(m)}_1+U^{(m)}_2+\ldots+U^{(m)}_j=l).
\end{align}
recalling the definition of $s_0$ and using the inequality $\log(1+x)\ge x-\frac{x^2}{2}$. $x\ge0$, we conclude that 
\begin{align}
\label{eq:B.2}
\nonumber
s_0^l
&=\exp\left\{l\log\left(1+\frac{2}{Bm}\left(1-\sqrt{\frac{j}{l}}\right)\right)\right\}\\
&\ge \exp\left\{\frac{2l}{Bm}\left(1-\sqrt{\frac{j}{l}}\right)
-\frac{2l}{B^2m^2}\right\}
\ge \frac{1}{C}\exp\left\{\frac{2l}{Bm}\left(1-\sqrt{\frac{j}{l}}\right)\right\}
\end{align}
provided that $l\le m^2$.

Next, recalling that 
$$
\frac{1}{1-a_m(s)}=\frac{1}{1-s}+\frac{Bm}{2},
$$
we have
$$
a_m(s)=1+\frac{s-1}{1-\frac{Bm}{2}(s-1)}.
$$
Consequently,
$$
a_m(s_0)
=1+\frac{\frac{2}{Bm}\left(1-\sqrt{\frac{j}{l}}\right)}
{1-\frac{Bm}{2}\frac{2}{Bm}\left(1-\sqrt{\frac{j}{l}}\right)}
=1+\frac{2}{Bm}\left(\sqrt{\frac{l}{j}}-1\right)
$$
and
$$
a_m^j(s_0)\le 
\exp\left\{\frac{2j}{Bm}\left(\sqrt{\frac{l}{j}}-1\right)\right\}.
$$
Combining this with \eqref{eq:B.2}, we conclude that
\begin{align}
\label{eq:B.3}
\nonumber
s_0^{-l}a_m^j(s_0)
&\le C\exp\left\{-\frac{2l}{Bm}\left(1-\sqrt{\frac{j}{l}}\right)
+\frac{2j}{Bm}\left(\sqrt{\frac{l}{j}}-1\right)\right\}\\
\nonumber
&=C\exp\left\{-\frac{2}{Bm}\left(\sqrt{l}-\sqrt{j}\right)^2\right\}
=C\exp\left\{-\frac{2(l-j)^2}{Bm\left(\sqrt{l}+\sqrt{j}\right)^2}\right\}
\\
&\le C\exp\left\{-\frac{(l-j)^2}{2Bml}\right\}
\end{align}
for all $j\le l\le m^2$. 

Thus, we are left to analyze the probability
$$
\pr(U^{(m)}_1+U^{(m)}_2+\ldots+U^{(m)}_j=l).
$$
It is easy to see that one has the following series representation
\begin{align}
\label{eq:am-gen}
a_m(s)=\frac{Bm}{2+Bm}+\frac{4}{Bm(2+Bm)}
\sum_{k=1}^\infty q_m^ks^k,
\end{align}
where 
$$
q_m=\frac{Bm}{2+Bm}.
$$
Therefore,
$$
\widehat{a}_m(s)=
\frac{Bm}{a_m(s_0)(2+Bm)}+\frac{4}{a_m(s_0)Bm(2+Bm)}
\sum_{k=1}^\infty q_m^ks_0^ks^k.
$$
This implies that, for every $k\ge 1$,
$$
\pr(U^{(m)}_1=k|U^{(m)}_1>0)
=\frac{4s_0^kq_{{m}}^k}{a_m(s_0)(1-\widehat{a}_m(0))Bm(2+Bm)}.
$$
Notice next that 
$$
s_0q_m
=\left(1+\frac{2}{Bm}\left(1-\sqrt{\frac{j}{l}}\right)\right)
\left(1-\frac{2}{Bm}\right)\le 1.
$$
Consequently,
\begin{align}
\label{eq:B.4}
\max_{k\ge 1}\pr(U^{(m)}_1=k|U^{(m)}_1>0)
\le \frac{C}{m^2(1-\widehat{a}_m(0))}.
\end{align}
By the definition of $\widehat{a}_m$,
\begin{align}
\label{eq:am-zero}
\nonumber
 1-\widehat{a}_m(0)
 &=1-\frac{a_m(0)}{a_m(s_0)}
 =\frac{(a_m(s_0)-1)+(1-a_m(0))}{a_m(s_0)}\\
 \nonumber
 &=\frac{\frac{2}{Bm}\left(\sqrt{\frac{l}{j}}-1\right)
 +\frac{2}{Bm+2}}{1+\frac{2}{Bm}\left(\sqrt{\frac{l}{j}}-1\right)}
 =\frac{\sqrt{\frac{l}{j}}-1+\frac{Bm}{Bm+2}}{\frac{Bm}{2}+\sqrt{\frac{l}{j}}-1}\\
 &\ge \frac{\sqrt{\frac{l}{j}}-\frac{2}{2+Bm}}{\frac{Bm}{2}+\sqrt{\frac{l}{j}}}
 \ge\frac{1}{C}\frac{\sqrt{\frac{l}{j}}}{m},
\end{align}
provided that $l\le m^2$.
Combining this with \eqref{eq:B.4}, we conclude that
\begin{align}
\label{eq:B.44}
\max_{k\ge 1}\pr(U^{(m)}_1=k|U^{(m)}_1>0)
\le \frac{C}{m},\quad j\le l
\end{align}
and that the right hand side is small for large $m$.

For each $m\ge1$, let $\{W_j^{(m)}\}_{j\ge1}$ be i.i.d. random variables with the distribution 
$$
\pr(W_1^{(m)}=k)=\pr(U^{(m)}_1=k|U^{(m)}_1>0),\quad k\ge1.
$$
By the total probability law,
\begin{align*}
&\pr(U^{(m)}_1+U^{(m)}_2+\ldots+U^{(m)}_j=l)\\
&\hspace{1cm}
=\sum_{t=1}^j \binom{j}{t}\widehat{a}_m^{j-t}(0)
(1-\widehat{a}_m(0))^j\pr(W_1^{(m)}+W_2^{(m)}+\ldots+W_t^{(m)}=l).
\end{align*}
Combining \eqref{eq:B.4} and \eqref{eq:B.44} with the classical bound on the concentration function of a sum of i.i.d. random variables, we have
\begin{align*}
\pr(W_1^{(m)}+W_2^{(m)}+\ldots+W_t^{(m)}=l)
&\le \frac{C_1}{\sqrt{t}}\frac{\max_{k\ge 1}\pr(W_1^{(m)}=k)}{1-\max_{k\ge 1}\pr(W_1^{(m)}=k)}\\
&\le\frac{C_2}{\sqrt{t}}\frac{1}{m^2(1-\widehat{a}_m(0))}.
\end{align*}
Consequently,
\begin{align*}
 &\pr(U^{(m)}_1+U^{(m)}_2+\ldots+U^{(m)}_j=l)\\
&\hspace{1cm}
\le \frac{C}{m^2(1-\widehat{a}_m(0))}
\sum_{t=1}^j\binom{j}{t}\widehat{a}_m^{j-t}(0)
(1-\widehat{a}_m(0))^j \frac{1}{\sqrt{t}}\\
&\hspace{1cm}
=\frac{C}{m^2(1-\widehat{a}_m(0))}
\e\left[\widehat{\zeta}_{j,m}^{-1/2};\zeta_{j,m}>0\right],
\end{align*}
where $\widehat{\zeta}_{j,m}$ has the binomial distribution with parameters $j$ and $(1-\widehat{a}_m(0))$. Applying now Lemma~18 from \cite{NV06}, we infer that
\begin{align}
\label{eq:B.5}
\pr(U^{(m)}_1+U^{(m)}_2+\ldots+U^{(m)}_j=l)
\le \frac{C}{m^2(1-\widehat{a}_m(0))}
\left(j(1-\widehat{a}_m(0))\right)^{-1/2}.
\end{align}
 Combining \eqref{eq:B.5} and \eqref{eq:am-zero}, we obtain
\begin{align}
\label{eq:B.6}
\pr(U^{(m)}_1+U^{(m)}_2+\ldots+U^{(m)}_j=l)
\le C\frac{j^{1/4}}{m^{1/2}l^{3/4}}.
\end{align}
Plugging \eqref{eq:B.3} and \eqref{eq:B.6} into \eqref{eq:B.1},
we obtain the desired bound.
\end{proof}
We next derive asymptotics for $\pr(G_m=k|G_0=j)$ in the case when the difference $k-j$ is not too big.
\begin{proposition}
\label{prop:lclt}
Assume that $k,m\to\infty$ and $\frac{k}{m}\to\infty$ sufficiently slow.
Then, for every fixed $A>0$,
$$
\pr(G_m=k|G_0=j)\sim\frac{1}{\sqrt{2\pi Bmk}}
\exp\left\{-\frac{(k-j)^2}{2Bmk}\right\}
$$
uniformly in $|j-k|\le A\sqrt{mk}$.
\end{proposition}
\begin{proof}
Set
$$
h=\frac{Bm}{2}.
$$
It is immediate from \eqref{eq:am-gen} that 
$$
\pr(G_m=0|G_0=1)=a_m(0)=\frac{Bm}{2+Bm}=(1+1/h)^{-1}
$$
and 
$$
\pr(G_m=k|G_0=1)=\frac{4}{Bm(2+Bm)}q_m^k
=\frac{1}{h^2}(1+1/h)^{-k-1},\quad k\ge1.
$$

For each fixed $m\ge 1$, let $\{V^{(m)}_j\}_{j\ge1}$ be independent, identically distributed random variables with the distribution given by
\begin{align}
\label{eq:A.1}
\pr(V^{(m)}_1=k)=\pr(G_m=k|G_0=1,G_m>0)
=\frac{1}{h}(1+1/h)^{-k},\quad k\ge1.
\end{align}

By the total probability law,
\begin{align*}
&\pr(G_m=k|G_0=j)\\
&\hspace{1cm}
=\sum_{t=1}^{k\wedge j}\binom{j}{t}(1-a_m(0))^ta_m^{j-t}(0)
\pr(V^{(m)}_1+V^{(m)}_2+\ldots+V^{(m)}_t=k)\\
&\hspace{1cm}
=(1+1/h)^{-j}\sum_{t=1}^{k\wedge j}\binom{j}{t} h^{-t}
\pr(V^{(m)}_1+V^{(m)}_2+\ldots+V^{(m)}_t=k).
\end{align*}
Next, \eqref{eq:A.1} implies that 
\begin{align*}
 \pr(V^{(m)}_1+V^{(m)}_2+\ldots+V^{(m)}_t=k)
 &=\sum_{k_j\ge 1:\sum k_j=k}\prod_{j=1}^t\pr(V_j^{(m)}=k_j)\\
 &=h^{-t}(1+1/h)^{-k}\sum_{k_j\ge 1:\sum k_j=k}1\\
 &=\binom{k-1}{t-1}h^{-t}(1+1/h)^{-k}.
\end{align*}
Consequently,
\[
\mathbb{P}(G_m = k \mid G_0 = j) = \left(1+\frac{1}{h}\right)^{-k-j} \sum_{t=1}^{k\wedge j} \binom{j}{t}\binom{k-1}{t-1} h^{-2t}.
\]
Let $f(x) = x\left(1+\frac{x}{h}\right)^{k-1}\left(1+\frac{1}{hx}\right)^{j}$, then
\[
[x^0] f(x) = \sum_{t=1}^{k\wedge j} \binom{j}{t}\binom{k-1}{t-1} \frac{1}{h^{2t-1}}.
\]
Hence,
\begin{align} \label{eq:A.2}
\mathbb{P}(G_m = k \mid G_0 = j) = \frac{1}{h}\left(1+\frac{1}{h}\right)^{-k-j} [x^0] f(x).
\end{align}
It holds
\[
[x^0] f = \frac{1}{2\pi}\int_{-\pi}^{\pi} e^{i\theta}\left(1+\frac{e^{i\theta}}{h}\right)^{k-1}\left(1+\frac{e^{-i\theta}}{h}\right)^{j} d\theta
\]
If $m\to\infty$ and $\frac{k}{m}\to\infty$ sufficiently slow then
\[
\left(1+\frac{e^{i\theta}}{h}\right)^{k-1}\left(1+\frac{e^{-i\theta}}{h}\right)^{j} \sim \exp\left\{e^{i\theta}\frac{k}{h} + e^{-i\theta}\frac{j}{h}\right\}
\]
uniformly in $\theta\in[-\pi,\pi]$. Then
\[
[x^0] f \sim \frac{1}{2\pi}\int_{-\pi}^{\pi} e^{i\theta} \exp\left\{e^{i\theta}\frac{k}{h} + e^{-i\theta}\frac{j}{h}\right\} d\theta.
\]
Since the series for the exponent converges absolutely, we have
\begin{align*}
\frac{1}{2\pi}\int_{-\pi}^{\pi} e^{i\theta} \exp\left\{a e^{i\theta} + b e^{-i\theta}\right\} d\theta 
&=
\frac{1}{2\pi}\sum_{p,q\ge 0}
    \frac{a^p b^q}{p!q!}
    \int_{-\pi}^\pi
        e^{i\theta} e^{i(p-q)\theta} d\theta\\
&=
\sum_{p=0}^\infty
\frac{a^p b^{p+1}}{p!(p+1)!}
=\sqrt{\frac{b}{a}}\, I_1\!\left(2\sqrt{ab}\right),
\end{align*}
where $I_1$ is modified Bessel function, see \cite[\href{https://dlmf.nist.gov/10.25.E2}{(10.25.2)}]{DLMF}. By \cite[\href{https://dlmf.nist.gov/10.40.E1}{(10.40.1)}]{DLMF}
\[
I_1(z) \sim \frac{e^{z}}{\sqrt{2\pi z}} \quad \text{as } |z|\to\infty \text{ and } |\arg z| \le \frac{\pi}{2}.
\]
For $a=k/h$ and $b=j/h$ we then obtain
\begin{align} \label{eq:A.3}
    [x^0]f \sim \sqrt{\frac{b}{a}}\, I_1\!\left(2\sqrt{ab}\right) \sim \sqrt{\frac{j}{k}}\, e^{2\frac{\sqrt{jk}}{h}} \frac{1}{\sqrt{2\pi}} \frac{\sqrt{h}}{k^{1/4} j^{1/4}} = \frac{j^{1/4} \sqrt{h}}{\sqrt{2\pi}\, k^{3/4}}\, e^{2\frac{\sqrt{jk}}{h}}
\end{align}
Notice that
\begin{align*}
\left(1+\frac{1}{h}\right)^{-k-j}
&=\exp\left\{-(k+j)\log\left(1+\frac{1}{h}\right)\right\}\\
&=\exp\left\{-\frac{(k+j)}{h}+O\left(\frac{k+j}{h^2}\right)\right\}
\sim \exp\left\{-\frac{(k+j)}{h}\right\}
\end{align*}
uniformly in $|j-k|\le A\sqrt{mk}$.
Using this and substituting \eqref{eq:A.3} into \eqref{eq:A.2}:
\[
\mathbb{P}(G_m = k \mid G_0 = j) \sim \frac{j^{1/4}}{\sqrt{2\pi}\sqrt{h}\cdot k^{3/4}}\, 
\exp\left\{
2\frac{\sqrt{jk}}{h}
-\frac{k+j}{h}
\right\}.
\]

To complete the proof it remains to notice that
\begin{align*}
\frac{k+j}{h}-2\frac{\sqrt{kj}}{h}
&=\frac{2(k-j)^2}{Bm(\sqrt{k}+\sqrt{j})^2}
=\frac{2(k-j)^2}{2Bmk}\left(1+\frac{1}{2}
\left(\frac{\sqrt{j}}{\sqrt{k}}-1\right)\right)^{-2}\\
&=\frac{2(k-j)^2}{2Bmk}
\left(1+O\left(\frac{|k-j|}{k}\right)\right)\\
&=\frac{2(k-j)^2}{2Bmk}+O\left(\frac{\sqrt{m}}{\sqrt{k}}\right)
=\frac{2(k-j)^2}{2Bmk}+o(1)
\end{align*}
uniformly in $|j-k|\le A\sqrt{mk}$.
\end{proof}

Combining Propositions~\ref{prop:comparison} and \ref{prop:lclt}, one obtains easily the following result.
\begin{corollary}
\label{cor:lclt}
If $k,m\to\infty$ and if $\frac{k}{m}\to\infty$ sufficiently slow then
$$
\pr(Z_m=k|Z_0=j)\sim\frac{1}{\sqrt{2\pi Bmk}}e^{-(k-j)^2/2Bmk}
$$
uniformly in $|j-k|\le A\sqrt{mk}$.
\end{corollary}

We next extend this local limit theorem to the process with immigration.
\begin{lemma}
\label{lem:lclt-immi}
If $k,m\to\infty$ and if $\frac{k}{m}\to\infty$ sufficiently slow then
$$
\pr(Y_m=k|Y_0=j)\sim\frac{1}{\sqrt{2\pi Bmk}}e^{-(k-j)^2/2Bmk}
$$
uniformly in $|j-k|\le A\sqrt{mk}$.
\end{lemma}
\begin{proof}
 Let us fix an arbitrary $\varepsilon>0$. By the total probability law,
 \begin{align}
 \label{eq:E.1}
 \nonumber 
 \pr(Y_m=k|Y_0=j)
 &=\sum_{l=0}^k\pr(Y_m=l|Y_0=0)\pr(Z_m=k-l|Z_0=j)\\
 \nonumber
 &=\sum_{l\le\varepsilon\sqrt{mk}}\pr(Y_m=l|Y_0=0)\pr(Z_m=k-l|Z_0=j)\\
 &\hspace{1cm}
 +\sum_{l>\varepsilon\sqrt{mk}}\pr(Y_m=l|Y_0=0)\pr(Z_m=k-l|Z_0=j).
 \end{align}
 The assumption $m=o(k)$ implies that $m=o(\sqrt{mk})$. Then, by Lemma~\ref{lem:loc-ld},
 $$
 \max_{l>\varepsilon\sqrt{mk}}\pr(Y_m=l|Y_0=0)
 =o\left(\frac{1}{\sqrt{mk}}\right)
 $$
 for every $\varepsilon>0$. This implies that
 \begin{align}
 \label{eq:E.2}  
 \sum_{l>\varepsilon\sqrt{mk}}\pr(Y_m=l|Y_0=0)\pr(Z_m=k-l|Z_0=j)
 =o\left(\frac{1}{\sqrt{mk}}\right). 
 \end{align}
 It follows from Corollary~\ref{cor:lclt} that there exists a depending on $A$ constant $C_A$ such that 
 \begin{align*}
 \sup_{|k-j|\le A\sqrt{mk},l\le\varepsilon\sqrt{mk}} 
 \left|\pr(Z_m=k-l|Z_0=j)-\pr(Z_m=k|Z_0=j)\right|
 \le C_A\frac{\varepsilon}{\sqrt{mk}}.
 \end{align*}
 Thus,
 \begin{align*}
 &\Bigg|
 \sum_{l\le\varepsilon\sqrt{mk}}\pr(Y_m=l|Y_0=0)\pr(Z_m=k-l|Z_0=j)\\
 &\hspace{2cm}- 
 \pr(Z_m=k|Z_0=j)\pr(Y_m\le\varepsilon\sqrt{mk})
 \Bigg|\le C_A\frac{\varepsilon}{\sqrt{mk}}.
 \end{align*}
 Noting that, due to the assumption $m=o(k)$, 
 $\pr(Y_m\le\varepsilon\sqrt{mk})\to1$ and using \eqref{eq:app.7}, we obtain
 \begin{align}
 \label{eq:E.3}
 \nonumber
 &\Bigg|
 \sum_{l\le\varepsilon\sqrt{mk}}\pr(Y_m=l|Y_0=0)\pr(Z_m=k-l|Z_0=j)
 -\pr(Z_m=k|Z_0=j)\Bigg|\\
 &\hspace{3cm}\le C_A\frac{\varepsilon}{\sqrt{mk}}
 +o\left(\frac{1}{\sqrt{mk}}\right).
 \end{align}
 Combining \eqref{eq:E.1} -- \eqref{eq:E.3}, we conclude that
 $$
 \Big|\pr(Y_m=k|Y_0=j)-\pr(Z_m=k|Z_0=j)\Big|
 \le C_A\frac{\varepsilon}{\sqrt{mk}}
 +o\left(\frac{1}{\sqrt{mk}}\right).
 $$
 Multiplying both sides by $\sqrt{mk}$ and letting $\varepsilon\to0$, we conclude that
 $$
 \Big|\sqrt{mk}\pr(Y_m=k|Y_0=j)-\sqrt{mk}\pr(Z_m=k|Z_0=j)\Big|=o(1)
 $$
 uniformly in $|j-k|\le A\sqrt{mk}$. Applying now Corollary~\ref{cor:lclt}, we obtain the desired relation.
\end{proof}
\section{Proof of the local limit theorem}
Let us first show that it suffices to prove the theorem for $Y_0=0$. By the total probability law,
\begin{align*}
\pr(Y_n=k)
&=\sum_{j=0}^\infty \pr(Y_0=j)\pr(Y_n=k|Y_0=j)\\
&\le \pr(Y_n=k|Y_0=0)+\sum_{j=1}^\infty \pr(Y_0=j)\max_{l\ge 1}\pr(Z_n=l|Z_0=j).
\end{align*}
By \eqref{eq:app.7},
\begin{align*}
 \sum_{j=n}^\infty \pr(Y_0=j)\max_{l\ge1}\pr(Z_n=l|Z_0=j)
 \le \frac{C}{n}\pr(Y_0\ge n)=o(1/n).
\end{align*}
Furthermore, by (45) from \cite{ShW},
\begin{align*}
 \sum_{j=1}^n \pr(Y_0=j)\max_{l\ge1}\pr(Z_n=l|Z_0=j)
 \le \frac{C}{n^2}\e[Y_0;Y_0<n]=o(1/n).
\end{align*}
Therefore,
\begin{align}
\label{eq:red.1}
\pr(Y_n=k)\le \pr(Y_n=k|Y_0=0)+o(1/n) 
\end{align}
uniformly in $k\ge 1$. 

Notice also that
\begin{align}
\label{eq:red.2}
\pr(Y_n=k)\ge \pr(Y_n=k|Y_0=0)\e\left[\left(f_n(0)\right)^{Y_0}\right]
=(1+o(1))\pr(Y_n=k|Y_0=0)
\end{align}
uniformly in $k$. Combining \eqref{eq:red.1} and \eqref{eq:red.2}, we infer that the case of general starting distribution will follow from the claim in the particular case $Y_0=0$.

From now on, we shall assume that $Y_0=0$.
To prove Theorem~\ref{thm:loc.lim} it is sufficient to show that 
\begin{align}
\label{eq:R.1}
\sup_{k\geq an }\left| n\pr\left(Y_{n}=k\right)-\frac{1}{\Gamma\left(\gamma\right)}\frac{2}{B}\left(\frac{2k}{Bn} \right)^{\gamma-1}e^{-2k/(Bn)} \right| \to 0
\end{align}
for every {\it fixed} $a>0$.
Noting next that 
$$
\sup_{x>0}x^\gamma e^{-x}=\left(\frac{\gamma}{e}\right)^\gamma
$$
and using the bound $\pr(Y_n=k)\le Ck^{-1}$, see \eqref{eq:Ubound-1}, we infer that 
$$
\sup_{k\geq Rn }\left| n\pr\left(Y_{n}=k\right)-\frac{1}{\Gamma\left(\gamma\right)}\frac{2}{B}\left(\frac{2k}{Bn} \right)^{\gamma-1}e^{-2k/(Bn)} \right|\le\frac{C}{R}.
$$
Consequently, for every $R>1$,
\begin{align*}
&\limsup_{n\to\infty} 
\sup_{k\geq an }\left| n\pr\left(Y_{n}=k\right)-\frac{1}{\Gamma\left(\gamma\right)}\frac{2}{B}\left(\frac{2k}{Bn} \right)^{\gamma-1}e^{-2k/(Bn)} \right|\\
&\hspace{1cm}
\le \limsup_{n\to\infty} 
\sup_{k\in[an,Rn] }\left| n\pr\left(Y_{n}=k\right)-\frac{1}{\Gamma\left(\gamma\right)}\frac{2}{B}\left(\frac{2k}{Bn} \right)^{\gamma-1}e^{-2k/(Bn)} \right|+\frac{C}{R}.
\end{align*}
This implies that \eqref{eq:R.1} will follow from
\begin{align}
\label{eq:R.2}
\sup_{k\in[an,Rn] }\left| n\pr\left(Y_{n}=k\right)-\frac{1}{\Gamma\left(\gamma\right)}\frac{2}{B}\left(\frac{2k}{Bn} \right)^{\gamma-1}e^{-2k/(Bn)} \right|\to0
\end{align}
for all $0<a<R$.

By the Markov property, for every $m\in\{1,2,\ldots,n\}$, 
\begin{equation}\label{eq2}
\pr(Y_{n}=k)=\sum_{j=0}^{\infty} \pr(Y_{n-m}=j)\pr(Y_{m}=k| Y_{0}=j).
\end{equation} 
Denote $N:=k-A\sqrt{km}$ and $M:=k+A\sqrt{km}$, where $A$ is a large  positive number.
We split \eqref{eq2} into three summands
$\Sigma_{1}$, $\Sigma_2$ and $\Sigma_3$
where $\Sigma_1$ is the sum over $j\in [0,N)$, $\Sigma_2$ is the sum over $j\in [N,M)$ and $\Sigma_3$ is the sum over $j\in [M,\infty)$. 
We shall also choose $m=m(n)$ such that $\frac{m}{n}\to0$ sufficiently slow.

To prove \eqref{eq:R.2} we shall analyze every term $\Sigma_i$ separately. We shall show that $\Sigma_1$ and $\Sigma_3$ are asymptotically negligible, and then determine the exact asymptotic behaviour of $\Sigma_2$. We start by bounding $\Sigma_3$.
\begin{proposition}\label{prop:sigma3}
If $m=o(k)$, then 
    \begin{equation*}
\limsup_{k \to \infty} k\Sigma_{3} \leq \frac{C}{A^{2}}.
    \end{equation*}
In particular, if $m=o(n)$ then, for every $a>0$, 
 \begin{equation*}
\limsup_{n\to \infty}\ n\sup_{k\ge an}\Sigma_{3} \leq \frac{C}{A^{2}a}.
    \end{equation*}
\end{proposition}
\begin{proof}
Applying \eqref{eq:Ubound-1} to $Y_{n-m}$, we obtain
\begin{equation}\label{eq:boundforsigma3}
\Sigma_3\leq \frac{C}{k}\sum_{j=M}^{\infty}\pr(Y_{m}=k | Y_{0}=j).
\end{equation}
Let $H_{m}(s)$ denote the generating function of $Y_{m}$ starting from $Y_{0}=0$. Then the generating function of the sequence 
$$
\theta_{k}:=\sum_{j=M}^{\infty}\pr(Y_{m}=k | Y_{0}=j),\quad k\ge0
$$
is given by 
$$
\sum_{j=M}^{\infty}f_{m}^{j}(s) H_{m}(s)=\frac{H_{m}(s)f^{M}(s)}{1-f_{m}(s)}.
$$
By Lemma 21 in \cite{NV06}, there exists a constant $C_{0}$ such that
$$
[s^k]\left[\frac{1}{1-f_{m}(s)}\right] \leq C_{0}
$$
uniformly in $k,m \geq 1$. This implies that 
\begin{align}\label{eq3}
\nonumber
\theta_{k}&=\sum_{j=0}^{k}[s^j]\left[ \frac{1}{1-f_{m}(s)}\right] 
[s^{k-j}]\left[H_{m}(s)f^{M}(s)\right]\\
\nonumber
&=\frac{1}{1-f_{m}(0)}\pr(Y_{m}=k | Y_{0}=M)+\sum_{j=1}^{k}[s^j]\left[\frac{1}{1-f_{m}(s)}\right]\pr(Y_m=k-j|Y_0=M)\\
&\leq \frac{1}{1-f_{m}(0)}\pr(Y_{m}= k | Y_{0}=M)+ C_{0}\pr(Y_{m}< k | Y_{0}=M).
\end{align}
Noting that $\pr(Y_{m} < k | Y_{0}=M)\leq \pr(Z_{m} < k | Z_{0}=M)$ and applying the Chebyshev inequality, we obtain
$$
\pr(Y_{m} < k | Y_{0}=M)\leq \frac{\operatorname{Var}[Z_{m}|Z_{0}=M]}{(M-k)^{2}}=\frac{MBm}{(M-k)^{2}}.
$$
By the choice of $M$, one has that
\begin{equation}\label{eq4}
    \pr(Y_{m}< k | Y_{0}=M) \leq B\frac{Mm}{A^{2}km}=\frac{B(1+A\sqrt{\frac{m}{k}})}{A^{2}}.
\end{equation}
To estimate the first summand on the right hand side of \eqref{eq3}, we notice that
\begin{align}\label{eq5}
\nonumber
\pr(Y_{m}=k| Y_{0}=M)
&=\sum_{l=0}^{k}\pr(Z_{m}=l| Z_{0}=M)\pr( Y_{m}=k-l |Y_{0}=0)\\
\nonumber
&=\pr(Z_{m}=0|Z_{0}=M) \pr(Y_{m}=k| Y_{0}=0)\\
&\hspace{1cm}
+\sum_{l=1}^{k}\pr(Z_{m}=l| Z_{0}=M)\pr(Y_{m}=k-l|Y_{0}=0).
\end{align}
Furthermore, due to \eqref{eq:Ubound-1},
\begin{equation}\label{eq6}
\pr(Z_{m}=0 | Z_{0}=M)\pr(Y_{m}=k| Y_{0}=0)\leq \frac{C}{k}(f_{m}(0))^{M} \leq \frac{C}{k}\exp\left\{-\gamma \frac{M}{m} \right\}
\end{equation}
and
$$
\sum_{l=1}^{k}\pr(Z_{m}=l|Z_{0}=M)\pr(Y_{m}=k-l|Y_{0}=0) \leq \sup_{l\geq 1}\pr(Z_{m}=l|Z_{0}=M).
$$

Plugging \eqref{eq:app.7} and \eqref{eq6} into \eqref{eq5}, we conclude that
$$
    \pr(Y_{m} = k | Y_{0}=M) \leq C \left( \frac{1}{k}\exp\left\{-\gamma\frac{M}{m}\right\}+\frac{1}{\sqrt{Mm}}\right).
$$
Applying this and \eqref{eq4} to the right-hand side of \eqref{eq3}, we finally obtain
\begin{equation}
\theta_{k}\leq \frac{C\left(1+A\sqrt{\frac{m}{k}}\right)} {A^{2}}+C\left(\frac{m}{k}\exp\left\{-\gamma \frac{k}{m}\right\}+\sqrt{\frac{m}{k}} \right).
\end{equation}
Hence, combining the latter inequality with \eqref{eq:boundforsigma3}, we arrive at the desired estimate.
\end{proof}
\begin{proposition}\label{prop:sigma1}
  Under the assumptions of Theorem 1, for every fixed $a>0$,
    \begin{equation*}
\limsup_{n\to\infty}n\sup_{k\ge an }\Sigma_1\le \frac{C}{A^2a}.
    \end{equation*}   
\end{proposition}
\begin{proof}
We first note that
\begin{align*}
\pr(Y_{m}=k| Y_{0}=j)
&=\sum_{l=0}^{k}\pr(Z_{m}=l| Z_{0}=j)\pr( Y_{m}=k-l |Y_{0}=0)\\
&\leq \max_{l\leq N_1} \pr(Y_{m}=k-l | Y_{0}=0) + 
\max_{l\in(N_1,k]} \pr(Z_{m}=l | Z_{0}=j),
\end{align*}
where
$$
N_1=N+\frac{A}{2}\sqrt{km}.
$$
Consequently,
\begin{align}
\label{eq:sigma1.1}
\nonumber
\Sigma_1
&=\sum_{j< N}\pr(Y_{n-m}=j)\pr(Y_{m}=k| Y_{0}=j)\\
&\le \max_{l\leq N_1} \pr(Y_{m}=k-l | Y_{0}=0)
+\sum_{j< N}\pr(Y_{n-m}=j)\max_{l\in(N_1,k]} \pr(Z_{m}=l | Z_{0}=j).
\end{align}
It follows from Lemma~\ref{lem:loc-ld} that 
\begin{align*}
\max_{l\leq N_1} \pr(Y_{m}=k-l | Y_{0}=0)
&=\max_{l\geq \frac{1}{2}A\sqrt{km}} \pr(Y_{m}=l | Y_{0}=0)\\
\le C
\frac{m+\displaystyle\max_{l\ge \frac{1}{2}A\sqrt{km}}\delta(l)}{A^2km}.
\end{align*}
Recalling that $\delta(l)=o(l)$, we infer that we may choose $m=m(n)=o(n)$ so that 
$$
\frac{\displaystyle\max_{k\ge an}\max_{l\ge \frac{1}{2}A\sqrt{km}}\delta(l)}{km}=o(1/n).
$$
This implies that
\begin{align}
\label{eq:sigma1.2}
\limsup_{n\to\infty}\ n\max_{k\ge an}\max_{l\leq N_1} \pr(Y_{m}=k-l | Y_{0}=0)\le \frac{C}{A^2a}.
\end{align}

To estimate the second part of the right hand side in \eqref{eq:sigma1.1}, we shall apply the approximation by geometric distributions. More precisely, it follows from Proposition~\ref{prop:comparison} that 
$$
\max_{l\in(N_1,k]} 
\left|\pr(Z_{m}=l | Z_{0}=j)-\pr(G_{m}=l | G_{0}=j)\right|
\le\frac{\beta(m)}{m^{3/2}}\sqrt{j}.
$$
Consequently,
\begin{align*}
&\sum_{j< N}\pr(Y_{n-m}=j)\max_{l\in(N_1,k]} \pr(Z_{m}=l | Z_{0}=j)\\
&\hspace{1cm}
\le \sum_{j< N}\pr(Y_{n-m}=j)\max_{l\in(N_1,k]} \pr(G_{m}=l | G_{0}=j)
+\frac{\beta(m)}{m^{3/2}}\e\sqrt{Y_n}.
\end{align*}
It is immediate from the Jensen inequality that $\e\sqrt{Y_n}\le C\sqrt{n}$. Since $\gamma(m)\to0$ as $m\to\infty$, we have
$$
\frac{\beta(m)}{m^{3/2}}\sqrt{n}=o(1/n),
$$
provided that $\frac{m}{n}\to0$ sufficiently slow. As a result we have
\begin{align}
\label{eq:sigma1.3}
\nonumber
&\sum_{j< N}\pr(Y_{n-m}=j)\max_{l\in(N_1,k]} \pr(Z_{m}=l | Z_{0}=j)\\
&\hspace{1cm}
\le \sum_{j< N}\pr(Y_{n-m}=j)\max_{l\in(N_1,k]} \pr(G_{m}=l | G_{0}=j)
+o(1/n).
\end{align}
It is immediate from Proposition~\ref{prop:geom-ub} that 
\begin{align}
\label{eq:B.7}
\max_{l\in[N_1,k]}\pr(G_m=l|G_0=j)
\le \frac{Cj^{1/4}}{m^{1/2}k^{3/4}}
\exp\left\{-\frac{(N_1-j)^2}{2Bmk}\right\}
\end{align}
for all $j<N$. This bound implies that, for all sufficiently large $k$,
\begin{align*}
&\sum_{j=1}^{k/2}\pr(Y_{n-m}=j)\max_{l\in[N_1,k]}\pr(G_m=l|G_0=j)\\
&\hspace{1cm}
\le \frac{C_1}{m^{1/2}k^{3/4}}
\exp\left\{-\frac{(N_1-k/2)^2}{2Bmk}\right\}
\sum_{j=1}^{k/2}j^{1/4}\pr(Y_{n-m}=j)\\
&\hspace{1cm}
\le \frac{C_1}{m^{1/2}k^{3/4}}
\exp\left\{-\frac{k}{9Bm}\right\}\e[Y_{n-m}^{1/4}].
\end{align*}
By the Jensen inequality,
$$
\e[Y_{n-m}^{1/4}]\le\left(\e Y_{n-m}\right)^{1/4}\le Cn^{1/4}.
$$
From this bound and from the assumption $k\ge an$ we obtain 
\begin{align}
\label{eq:B.8}
\nonumber
&\sum_{j=1}^{k/2}\pr(Y_{n-m}=j)\max_{l\in[N_1,k]}\pr(G_m=l|G_0=j)\\
&\hspace{2cm}
\le \frac{C}{a^{3/4}\sqrt{mn}}\exp\left\{-\frac{an}{9Bm}\right\}
=o(n^{-1}).
\end{align}
Using \eqref{eq:B.7} once again and applying the bound
$\pr(Y_{n-m}=j)\le C/j$, we obtain 
\begin{align*}
&\sum_{j=k/2}^N\pr(Y_{n-m}=j)\max_{l\in[N_1,k]}\pr(G_m=l|G_0=j)\\
&\hspace{1cm}
\le \frac{C}{m^{1/2}k^{3/4}}\sum_{j=k/2}^N 
j^{-3/4}\exp\left\{-\frac{(N_1-j)^2}{2Bmk}\right\}\\
&\hspace{1cm}
\le \frac{C}{m^{1/2}k^{3/2}}
\sum_{j=N_1-N}^\infty \exp\left\{-\frac{j^2}{2Bmk}\right\}.
\end{align*}
Bounding the sum by the integral, one gets easily
$$
\sum_{j=N_1-N}^\infty \exp\left\{-\frac{j^2}{2Bmk}\right\}
\le \int_{(N_1-N)/2}^\infty e^{-x^2/2Bmk}dx.
$$
Recalling that $N_1-N=\frac{A}{2}\sqrt{mk}$ and using the standard bound for the gaussian integral, we conclude that 
$$
\sum_{j=N_1-N}^\infty \exp\left\{-\frac{j^2}{2Bmk}\right\}
\le \frac{C\sqrt{mk}}{\sqrt{A}}e^{-A^2/8B}.
$$
Consequently,
\begin{align*}
\sum_{j=k/2}^N\pr(Y_{n-m}=j)\max_{l\in[N_1,k]}\pr(G_m=l|G_0=j)
\le\frac{C}{k\sqrt{A}}e^{-A^2/8B}.
\end{align*}
Combining this with \eqref{eq:B.8} and recalling that $k\ge an$,
we conclude that
$$
\limsup_{n\to\infty}\ n\sum_{j< N}\pr(Y_{n-m}=j)\max_{l\in(N_1,k]} \pr(G_{m}=l | G_{0}=j)
\le
\frac{C}{a\sqrt{A}}e^{-A^2/8B}.
$$
Combining this with \eqref{eq:sigma1.3} finished the proof.
\end{proof}

We finally turn to the analysis of $\Sigma_2$.
\begin{proposition}
\label{prop:sigma2}
For all $a,A,R>0$, uniformly in $k\in[an, Rn]$,
$$
\lim_{n\to\infty}\frac{n}{p(k/n)}\Sigma_2=
\frac{1}{\sqrt{2\pi}}\int_{-A}^{A}e^{-u^2/2}du,
$$
where 
$$
p(x)=\frac{1}{\Gamma(\gamma)}\left(\frac{2}{B}\right)^\gamma x^{\gamma-1}e^{-2x/B},\ x>0
$$
and
$$
\gamma=\frac{2\lambda}{B}.
$$
\end{proposition}
\begin{proof}
According to the limit theorem for $Y_n/n$, there exists a sequence $\delta_n\downarrow0$ such that
$$
\sup_{x\ge 0}
\left|\pr\left(\frac {Y_n}{n}\le x\right)-F(x)\right|
=o(\delta_n),
$$
where $F$ is the distribution function of the distribution with the density $p(x)$. This implies that 
$$
\sup_{y\ge x\ge 0}
\left|\pr\left(x\le \frac{Y_n}{n}< y\right)-\int_x^y p(u)du\right|
=o(\delta_n).
$$
Fix some $x_0>0$. Since $p(x)$ is continuous and bounded on $[x_0,\infty)$, we conclude that, for any sequence $\theta_n\to0$,
\begin{align}
\label{eq:rate_of_conv}
\sup_{x\ge x_0, \delta_n\le y-x\le \theta_n}
\left|\frac{\pr\left(x\le \frac{Y_n}{n}< y\right)}{y-x}-p(x)\right|\to0.
\end{align}
Fix some $\varepsilon>0$. We may assume, without loss of generality, that 
$$
\frac{\sqrt{mn}}{n\delta_{n/2}}\to\infty.
$$
(If this relation is not valid, we can simply increase $m$.) This allows us to apply \eqref{eq:rate_of_conv} and to obtain 
\begin{align*}
\pr\left(Y_{n-m}\in
\left[k+l\varepsilon\sqrt{Bmk},k+(l+1)\varepsilon\sqrt{Bmk})\right)
\right) 
&\sim \frac{\varepsilon\sqrt{Bmk}}{n-m}
p\left(\frac{k+l\varepsilon\sqrt{Bmk}}{n-m}\right)\\
&\sim  \frac{\varepsilon\sqrt{Bmk}}{n}p\left(\frac{k}{n}\right)
\end{align*}
for every fixed $l$.

Noting that 
$$
e^{-(l+1)^2\varepsilon^2/2}
\le e^{-(k-j)^2/2Bmk}\le e^{-l^2\varepsilon^2/2}
$$
for all $j\in\left[k+l\varepsilon\sqrt{Bmk},k+(l+1)\varepsilon\sqrt{Bmk})\right)$, we conclude that
\begin{align*}
&\sum_{j=k}^M\pr(Y_{n-m}=j)\frac{1}{\sqrt{2\pi Bmk}}e^{-(k-j)^2/2Bmk}\\
&\hspace{5mm}
\le \frac{\varepsilon}{\sqrt{2\pi Bmk}}\sum_{l\le A/\varepsilon}
e^{-l^2\varepsilon^2/2}\pr\left(Y_{n-m}\in
\left[k+l\varepsilon\sqrt{Bmk},k+(l+1)\varepsilon\sqrt{Bmk})\right)
\right)\\
&\hspace{5mm}
=\frac{(1+o(1))p(k/n)}{n}\frac{\varepsilon}{\sqrt{2\pi}}
\sum_{l\le A/\varepsilon}e^{-l^2\varepsilon^2/2}
\end{align*}
and
\begin{align*}
&\sum_{j=k}^M\pr(Y_{n-m}=j)\frac{1}{\sqrt{2\pi Bmk}}e^{-(k-j)^2/2Bmk}\\
&\hspace{5mm}
\ge \frac{\varepsilon}{\sqrt{2\pi Bmk}}\sum_{l\le A/\varepsilon-1}
e^{-(l+1)^2\varepsilon^2/2}\pr\left(Y_{n-m}\in
\left[k+l\varepsilon\sqrt{Bmk},k+(l+1)\varepsilon\sqrt{Bmk})\right)
\right)\\
&\hspace{5mm}
=\frac{(1+o(1))p(k/n)}{n}\frac{\varepsilon}{\sqrt{2\pi}}
\sum_{l\le A/\varepsilon-1}e^{-(l+1)^2\varepsilon^2/2}.
\end{align*}
It is easy to see that, as $\varepsilon\to0$, 
$$
\frac{\varepsilon}{\sqrt{2\pi}}
\sum_{l\le A/\varepsilon-1}e^{-(l+1)^2\varepsilon^2/2}\to 
\frac{1}{\sqrt{2\pi}}\int_0^A e^{-u^2/2}du
$$
and 
$$
\frac{\varepsilon}{\sqrt{2\pi}}
\sum_{l\le A/\varepsilon}e^{-l^2\varepsilon^2/2}\to 
\frac{1}{\sqrt{2\pi}}\int_0^A e^{-u^2/2}du.
$$
Consequently,
\begin{align*}
\sum_{j=k}^M\pr(Y_{n-m}=j)\frac{1}{\sqrt{2\pi Bmk}}e^{-(k-j)^2/2Bmk}=\frac{(1+o(1))p(k/n)}{n}\frac{1}{\sqrt{2\pi}}\int_0^A e^{-u^2/2}du.
\end{align*}
By the same argument,
\begin{align*}
\sum_{j=N}^{k-1}\pr(Y_{n-m}=j)\frac{1}{\sqrt{2\pi Bmk}}e^{-(k-j)^2/2Bmk}
=\frac{(1+o(1))p(k/n)}{n}\frac{1}{\sqrt{2\pi}}\int_{-A}^0 e^{-u^2/2}du.
\end{align*}
Combining this relations with Lemma~\ref{lem:lclt-immi}, we have
\begin{align*}
\Sigma_2
&\sim\sum_{j=N}^M \pr(Y_{n-m}=j)\frac{1}{\sqrt{2\pi Bmk}}e^{-(k-j)^2/2Bmk}\\
&\sim \frac{p(k/n)}{n}\frac{1}{\sqrt{2\pi}}\int_{-A}^A e^{-u^2/2}du.
\end{align*}
This completes the proof of the proposition.
\end{proof}

We are now in position to complete the proof of the local limit theorem. We have
\begin{align*}
    \left|
    n \pr(Y_n = k) - p(k/n)
    \right|
\le
    n\Sigma_1
+
    n\Sigma_3
&+
    \left|
        n\Sigma_2 
    -
        \frac{p(k/n)}{\sqrt{2\pi}}
        \int_{-A}^A e^{-u^2/2}du    
    \right|\\
&\hspace{3cm}+
    \sqrt{\frac2\pi}
    \int_A^{\infty}
        e^{-u^2/2}du.
\end{align*}
Taking supremum over $k \in [an,Rn]$, considering $\limsup$ and applying Propositions \ref{prop:sigma1}, \ref{prop:sigma2}, \ref{prop:sigma3}, we obtain
\begin{align*}
\limsup_{n\to \infty}
\sup_{ k \in [an, Rn]}
\left|
    n \pr(Y_n = k) - p(k/n)
    \right|
\le
    \frac{C}{A^2a}
+
    \frac{2}{A}e^{-A/2}.
\end{align*}
Since $A$ is arbitrary, the limit on the right-hand side is zero. This gives us \eqref{eq:R.2}. Thus, Theorem~\ref{thm:loc.lim} is proved.

\end{document}